\documentclass[reqno]{amsart}
\usepackage{latexsym}
\usepackage{amssymb,amsthm,amsmath}

\usepackage{color}
\usepackage[
breaklinks,colorlinks=true]{hyperref}
\hypersetup{urlcolor=blue, citecolor=blue}
\usepackage{tikz}
\usepackage{multicol}

\def\dsp{\displaystyle}

\theoremstyle{plain}
\newtheorem{theorem}
{Theorem}
\newtheorem{lemma}
{Lemma}
{Corollary}
\newtheorem{proposition}
{Proposition}
{Example}
\newtheorem{definition}
{Definition}

\theoremstyle{definition}
\theoremstyle{remark}
\newtheorem{remark}
{Remark}

\def\d#1{{#1\kern-0.4em\char"16\kern-0.1em}}
\def\D#1{{\raise0.2ex\hbox{-}\kern-0.4em #1}}
\newcounter{zd}

\newcounter{zdr}[subsection]

\newcommand{\eps}{\varepsilon}
\def\mx{{\bf x}}

\def\my{{\bf y}}

\def\R{I\!\!R}

\def\pa{\partial}

\def\cal{\mathcal}

\DeclareMathOperator{\Div }{div}

\def\R{I\!\!R}

\def\cal{\mathcal}

\def\pa{\partial}

\def\b{\infty}

\def\be{\begin{equation}}
\def\ee{\end{equation}}

\def\mff{{\mathfrak f}}

\begin{document}

\title[Scalar conservation law with discontinuous flux revisited]
{
Entropy conditions for scalar conservation
laws\\ with discontinuous flux revisited}

\author{B. Andreianov}\address{Boris Andreianov,
Laboratoire de Math\'ematiques CNRS UMR 6623, Universit\'e de Franche-Comt\'e,
25030 Besan\c{c}on Cedex, France \;\;\;
and\;\;\; Institut f\"ur Mathematik, Technische Universit\"at Berlin,
Stra\ss e des 17. Juni 136, 10623 Berlin, Germany}
\email{boris.andreianov@univ-fcomte.fr}

\author{D.~Mitrovi\'c}\address{Darko Mitrovi\'c,
University of Montenegro, Faculty of Mathematics, Cetinjski put bb,
81000 Podgorica, Montenegro}\email{matematika@t-com.me}

\date{\today}



\begin{abstract}
We propose new entropy admissibility conditions for multidimensional
hyperbolic scalar conservation laws with discontinuous flux which generalize
 one-dimensional Karlsen-Risebro-Towers entropy conditions.
 These new conditions are designed, in particular, in order to characterize the limit of vanishing viscosity approximations.
 On the one hand, they  comply quite naturally with a certain class of physical and numerical modeling assumptions;  on the other hand, their mathematical assessment turns out to be intricate.

\smallskip
The generalization we propose is not only with respect to the space dimension, but mainly in the sense that the
``crossing condition'' of [K.H. Karlsen, N.H. Risebro, J. Towers, Skr.\,K.\,Nor.\,Vid.\,Selsk. (2003)] is not mandatory for proving uniqueness with the new definition.
We  prove uniqueness of solutions and give tools to justify their existence via the vanishing viscosity method,
for the multi-dimensional spatially inhomogeneous case with a finite number of Lipschitz regular hypersurfaces
 of discontinuity for the flux function.
\end{abstract}

\maketitle

\tableofcontents

\section{Introduction}

Conservation laws of the form
\begin{equation}\label{eq:cons-law}
\partial_t u + \text{div}_\mx \mff (t,\mx,u) = S(t,\mx,u)
\end{equation}
serve as mathematical models for one-dimensional gas dynamics, road
traffic, for flows in porous media with neglected capillarity
effects, blood flow, radar shape-from-shading problems, and in
several other applications. The multi-dimensional conservation law
also appears in coupled models, although in this case the regularity
of the flux function $\mff$ in $(t,\mx)$ is often not sufficient to
develop a full well-posedness theory. The mathematical theory of
\eqref{eq:cons-law} is very delicate because, in general, even for
regular data classical solutions need not to exist; on the other
hand, weak (distributional) solutions are, in general, not unique.
The classical theory is best established for Cauchy and
boundary-value problems in the case where $\mff$ is  Lipschitz
continuous in $(t,\mx)$ and uniformly locally Lipschitz continuous in $u$.
 The source $S$ can be assumed, e.g., Lipschitz continuous in $u$ uniformly in $(t,\mx)$. In
this case the S.N. Kruzhkov definition of entropy solution in
\cite{Kru} and the associated analysis techniques (vanishing
viscosity approximation for the existence proof, doubling of
variables for the uniqueness proof) provide a well-posedness
framework for \eqref{eq:cons-law}.

\subsection{Discontinuous-flux models and rough entropy inequalities}\label{ssec:INTRO-rough}

Local Lipschitz assumption on $k\mapsto \mff(t,\mx,k)$ is natural in many applications,
but the assumption of regular dependence on the spatial variable $\mx$ is very restrictive. Indeed, road traffic
with variable number of lanes (\cite{BKT2009-multilanes}), Buckley-Leverett equation in a layered porous medium (\cite{Kaas,AndrCances}), sedimentation applications (\cite{Die1,Die2,BurgerEtAl})
make appear models with piecewise regular, jump discontinuous in $\mx$ flux functions.
The theory of such problems, called \emph{discontinuous-flux conservation laws},
has been an intense subject of research in the last twenty years.
  The main goal of this research was to design
a suitable approach to definition and numerical approximation of entropy solutions, in relation with the physical context of
different discontinuous-flux models. Speaking of ``notion of solution'' in this context, one usually means
weak solutions subject to some additional \emph{admissibility conditions}, cf.~\cite{Kru}.

\smallskip
Almost all  admissibility conditions designed in the literature were confined to the one-dimensional case. We mention here the minimal jump condition \cite{GR}, minimal variation condition
and $\Gamma$-condition \cite{Die1,Die2}, entropy conditions
\cite{kar3,AG}, vanishing capillarity limit \cite{Kaas,AndrCances},
admissibility conditions via adapted entropies \cite{AP,kar4} or
via admissible jumps description at the interface \cite{sid_2,sid,Die3}. An
extensive overview on the subject as well as a kind of unification
of the mentioned approaches in the one-dimensional case was given in \cite{AKR}, where further references can be found.

\smallskip
Adimurthi et al.~\cite{sid_2} observed that infinitely many different, though equally
mathematically consistent, notions of solution may co-exist in the
discontinuous-flux problems; therefore the choice of solution notion
is a part of modeling procedure (see, e.g., \cite{AndrCances} for an
exhaustive study of the vanishing capillarity limits of the
one-dimensional Buckley-Leverett equation, where different sets of admissibility conditions
 are put forward for different choices of
physically relevant vanishing capillarity).
In the present contribution we limit our attention to
characterization of vanishing viscosity limit solutions for problems
of kind \eqref{eq:cons-law}; these approximations were studied in a huge number of works
(see, e.g., \cite{GR,Die1,Die2,towers1,towers2,SeguinVovelle,kar3,ken_chi,BachmannVovelle,kar4,pan_arma,Die3}) including several works in multiple space dimensions (\cite{JimenezLevi,Jimenez,AKRnhm,Pan_08,prag}) and they
remain relevant in several models based on discontinuous-flux conservation laws.

\smallskip
The basis of the different definitions of admissibility of solutions is provided by Kruzhkov entropy inequalities
(\cite{Kru}) in the regions of smoothness of the flux; the main difficulty consists in taking into account the jump discontinuities  of the flux.  To do so, for a quite general setting one may only assume that
 \begin{equation}\label{eq:BVflux}
 \text{for all $k\in\R$, $\mff(.,.,k)\in BV_{loc}(\R^+\times \R^d)$.}
 \end{equation}
This rather weak regularity appears naturally e.g. in the study of
triangular systems of conservation laws (see \cite{RisebroEtAl} and
references therein). In the framework \eqref{eq:BVflux}, under a
non-degeneracy assumption of the fluxes $k\mapsto \mff(t,\mx,k)$,
existence of solutions satisfying the family of entropy inequalities
  \begin{multline}\label{eq:Panov-def}
    \forall k\in\R\;\;\; |u\!-\!k|_t \; +\; \text{div}_\mx\Bigl(\text{sgn}(u\!-\!k)(\mff(t,\mx,u)\!-\!\mff(t,\mx,k)) \Bigr)
     - \text{sgn}(u\!-\!k)S(t,\mx,u) \\
    \leq\; - \text{sgn}(u\!-\!k)(\text{div}_\mx \mff(t,\mx,k))^{ac}\;+\;
    |(\text{div}_\mx \mff(t,\mx,k))^s|\;\; \text{in $\mathcal D'(\R^+_*\times\R^d)$}
  \end{multline}
  has been proved by Panov in \cite{pan_arma} using vanishing viscosity method; here,
  $$
 \text{div}_\mx \mff(t,\mx,k)=(\text{div}_\mx \mff(t,\mx,k))^{ac}+ (\text{div}_\mx \mff(t,\mx,k))^s
  $$
  is the Jordan decomposition of the Radon measure $\text{div}_\mx \mff(t,\mx,k)$ into its absolutely continuous part and its singular part
(cf. \eqref{sep827} below, for a particular but representative
case). Let us stress that inequalities \eqref{eq:Panov-def} use a
roughly estimated contribution of the jump singularities in the flux
$\mff$, which turns out to be a serious obstacle for proving
uniqueness of solutions in the sense~\eqref{eq:Panov-def}.

\smallskip
In order to explain this more accurately, let us consider the
one-dimensional variant of \eqref{eq:cons-law} augmented with the
initial conditions $u|_{t=0}=u_0 \in L^\infty(\R^d)$.
Condition \eqref{eq:Panov-def} is the general form of the
Karlsen-Risebro-Towers admissibility condition (see \cite{towers1,towers2,kar3}; see also~\cite{SeguinVovelle,BachmannVovelle}):
\begin{multline}\label{eq:KRT-cond}
\forall k\in [a,b]\;\; \pa_t {|u\!-\!k|} +\pa_x \Big\{{\rm sgn}(u-k)\Big[H(x)({f}(u)\!-\!{f}(k))+H(-x)({g}(u)\!-\!{g}(k))\Big]\Big\} \\
-|{f}({k})\!-\!{g}({k})|\delta_0(x) \leq 0 \
\ {\rm in} \ \ {\cal D}'(\R^+\times \R).
\end{multline}
This is precisely \eqref{eq:Panov-def} written for the case with source $S\equiv 0$ and flux
$\mff(t,x;u)=f(u)H(x)+g(u)H(-x)$; here and throughout the paper,
$$
\begin{array}{l}
\text{$H(\cdot)=\text{sgn}^+(\cdot)$ is the Heavyside function;}\\
 \text{$\delta_0$ is the Dirac measure concentrated at $0$, i.e., $\delta_0=H'$ in $\mathcal D'(\R)$.}
 \end{array}
$$
The advantage of definition of solution by the family of inequalities \eqref{eq:KRT-cond}
 lies in a simple and natural existence and stability proof. Using the vanishing viscosity approximations of
\eqref{eq:cons-law} for which existence is easy even in the
discontinuous-flux setting,
one readily obtains for the corresponding solutions $u_\eps$ the inequalities
\begin{multline}\label{eq:entropy-eps}
  \pa_t {|u_\eps\!-\!k|}+\pa_x \Big\{{\rm sgn}(u_\eps\!-\!k)\Big[H(x)({f}(u_\eps)\!-\!{f}(k))+H(-x)({g}(u_\eps)\!-\!{g}(k))\Big]\Big\}\\
  \dsp \qquad-\text{sgn}(u_\eps(t,0)\!-\!k)({f}({k})\!-\!{g}({k}))\delta_0(x) \leq \eps \pa_{xx} |u_\eps\!-\!k| \
\ {\rm in} \ \ {\cal D}'(\R^+\times \R),
\end{multline}
a.e. on $(0,T)$.
Roughly estimating the contribution of the term concentrated on $\{x=0\}$, we can pass to the limit
and get  \eqref{eq:KRT-cond}. In a similar way, stability of solutions in the sense \eqref{eq:KRT-cond} under a.e. convergence can be justified.

\smallskip
Concerning the question of uniqueness of solutions admissible in the Karlsen-Risebro-Towers sense, the following facts are known:

\smallskip
\noindent$\bullet$ Given a datum $u_0$, uniqueness of an admissible solution in the
sense \eqref{eq:KRT-cond} to the above problem holds true, provided $f,g$ satisfy the \emph{crossing condition}
\begin{equation}\label{eq:crossing-cond}
f(u)-g(u)<0<f(v)-g(v) \;\;\Rightarrow\;\; u<v
\end{equation}
on the shape of the fluxes $f$ and $g$. For the proof,
 see Towers \cite{towers1,towers2} (see also \cite{SeguinVovelle})
  for the case where the fluxes $f,g$ have no crossing points in $(a,b)$;
   see Karlsen, Risebro and Towers \cite{kar3} for general uniqueness argument under assumption
\eqref{eq:crossing-cond}.

\smallskip
\noindent$\bullet$
 If the crossing condition \eqref{eq:crossing-cond} fails, then at least for some initial data there exist more than one admissible solution in the sense  \eqref{eq:KRT-cond}, see \cite[Sect.\,4.7]{AKR}.



\smallskip
Several attempts were made already to improve the Karlsen-Risebro-Towers conditions.
Let us mention here the following directions.

\smallskip
\noindent 1. In the work \cite{dm_dcds}, the idea was to use a transformation of the original
equation (cf. \cite{Pan_08} and \cite{prag}) in order to enforce the crossing condition for the fluxes depending on the transformed unknown; yet such transformation is often artificial with respect to the underlying model.

\smallskip
\noindent  2. In \cite{Die3,AKR,AKRnhm} the crossing condition was dropped. \\[1pt]
(2a) In \cite{Die3}, a new version of the $\Gamma$-condition on the solution jumps was proposed, under which the uniqueness proof was achieved without the crossing condition. Admissible solutions
were characterized in terms of generalized Oleinik jump inequalities.\\[1pt]
(2b)  A different though equivalent to \cite{Die3} admissibility condition was proposed in \cite{AKR,AKRnhm} in terms of the \emph{vanishing viscosity germ} $\mathcal G_{VV}$ recalled in \S~\ref{ssec:FLAT-GVV-uniq} below.
It was shown that $\mathcal G_{VV}$-entropy solutions are always unique, regardless the shape of the fluxes.  The way to express this admissibility condition is rather tricky. For a straightforward approach put forward in~\cite{AKR}, one uses intricate interface coupling conditions. Alternatively, a carefully selected family of \emph{adapted entropies} (see \cite{BaitiJenssen,AP} and \cite{kar4}) can be used to characterize admissible solutions, in the place of the classical Kruzhkov entropies exploited in \eqref{eq:KRT-cond}. Neither the coupling conditions encoded in the germ $\mathcal G_{VV}$, nor the associated choice of adapted entropies are self-evident; their relation to the vanishing viscosity approximation follows from a lengthy analysis of possible viscosity profiles.

\smallskip
To sum up, none of the aforementioned approaches of admissibility is as intuitive and appealing as \eqref{eq:Panov-def} or \eqref{eq:KRT-cond}.
In the present paper, our goal is to give a definition of solution to \eqref{eq:cons-law} which, similarly to the definitions \eqref{eq:Panov-def} or \eqref{eq:KRT-cond}, could be seen as a natural one and that would lead to well-posedness without assuming the crossing condition.
Notice that in passing, we give an equivalent and somewhat more natural characterization of $\mathcal G_{VV}$-entropy solutions of \cite{AKR,AKRnhm}.

\subsection{Singular values of $u$ and refined entropy inequalities}\label{ssec:INTRO-refined}
The idea of this paper is to strengthen the definitions of
\cite{kar3,pan_arma} by suggesting a finer way to take into account the
contribution of the flux jumps into the entropy inequality. To
explain the idea, we go back to the general framework
\eqref{eq:BVflux}; in this case, we suggest to replace
\eqref{eq:Panov-def} by the following less restrictive inequalities:
\begin{multline}\label{eq:Panov-def-extended}
    \forall k\in\R\;\; |u\!-\!k|_t \; +\; \text{div}_\mx\Bigl(\text{sgn}(u\!-\!k)(\mff(t,\mx,u)\!-\!\mff(t,\mx,k)) \Bigr)\\
    \leq\; - \text{sgn}(p_u\!-\!k)\text{div}_\mx \mff(t,\mx,k)\;\; \text{in $\mathcal D'(\R^+\times\R^d)$}
  \end{multline}
where $p_u$ is some globally defined Borel function satisfying the property
\begin{equation}\label{eq:p-function}
 p_u=u \;\; \text{a.e. with respect to the Lebesgue measure}.
\end{equation}
Introducing everywhere defined $p_u$ means that we tacitly assign a
value to $u$ also on the singular set\footnote{The literature on
discontinuous-flux problems concentrates on the case where the union
in $k\in\R$ of the sets of singularities of $(t,\mx)\mapsto \text{div}_\mx
\mff(t,\mx,k)$ is of Lebesgue measure zero, typically it is a locally finite union $\Sigma$ of hypersurfaces in $\R^+\times\R^d$. In this case
 one can say that $u$ is defined on the set of full measure excluding singularities, and $p_u$ extends
 $u$ to the singular set $\Sigma$.} of $(t,\mx)$ where the flux $(t,\mx)\mapsto \mff(t,\mx,k)$
 experiences jumps or Cantor-type singularities. In the sequel, we will say that $p_u$ encodes the \emph{singular values of $u$}.

\begin{remark}\label{rem:micromodel} Let us present several observations concerning $p_u$.\\
(i) Introduction of unknown singular values $p_u$ on the jump manifolds $\Sigma$ may appear as a naive
 attempt to resolve the difficulty of definition of solution: indeed, in general these values cannot be observable.
 %

 \smallskip\noindent
 (ii)
 Yet from the modeling viewpoint, existence of singular values $p_u$ on $\Sigma$ can be put in relation with the fact that the vanishing viscosity approximation ensures global continuity of the approximate solutions $u_\eps$.
This approximation is suitable for models where
 the solution
  is expected to be continuous across $\Sigma$ at a finer spatial scale, undergoing rapid transition of interface layer kind.
Such fine-scale continuity assumptions are also natural in numerical approximation and modeling: we refer to \cite{ChavantCohenJaffre} for an early work based on this idea, and to
 \cite{AndrCances-transmission} for a deeper discussion of a wide class of related modeling hypotheses.

%
 \smallskip\noindent
 (iii) In practice, whenever the values of approximate solutions $u_\eps$
on a jump manifold $\Sigma$ of the flux $(t,\mx)\mapsto \mff(t,\mx,k)$
happen to exhibit a non-oscillatory behavior, one obtains $p_u$ on
$\Sigma$ as the pointwise limit of $u_\eps|_\Sigma$, for some vanishing sequence of $\eps$.
Mathematically, \emph{a priori} justification of strong convergence
of $u_\eps|_\Sigma$ to some limit $p_u$ seems out of reach even when
strong compactness of $(u_\eps)_\eps$ in $L^1_{loc}$ topology is
easy to justify.

\smallskip\noindent
(iv) It is also clear that singular values $p_u$ need
not be uniquely defined with respect to the Hausdorff measure on
$\Sigma$ in its natural dimension $d$: to observe this, in is enough
to consider the simplest case of converging one-dimensional
viscosity profiles that may provide a continuum of different values for
$p_u$. Proposition~\ref{ekv} suggests some canonical choice of $p_u$.
\end{remark}

\subsection{A brief description of technical ideas and obtained results}\label{ssec:INTRO-results}

Throughout the paper we assume that
\begin{equation}\label{eq:time-orthog-Sigma}
\begin{array}{l}
\text{at a.e. point of $\Sigma$ (with respect to the $d$-dimensional Hausdorff measure),}\\
 \text{$\Sigma$ is not orthogonal to the direction of the time axis.}
\end{array}
\end{equation}

As a matter of fact, we are not able to suggest technical
tools that would permit to exploit definition
\eqref{eq:Panov-def-extended} in the general setting
\eqref{eq:BVflux},\eqref{eq:time-orthog-Sigma} neither in view of existence nor in view of
uniqueness of solutions.
The results of this note only concern the practical case of
piecewise Lipschitz, jump discontinuous across a union  $\Sigma$ of Lipschitz hypersurfaces, and genuinely nonlinear fluxes.

\smallskip
In this case, firstly, we are able to prove uniqueness of admissible solutions
in the sense \eqref{eq:Panov-def-extended},\eqref{eq:p-function}
without any artificial condition on the flux crossing.
Uniqueness result for solutions in the sense
of~\eqref{eq:Panov-def-extended},\eqref{eq:p-function} is based upon consideration of
constraints that are imposed by inequalities
\eqref{eq:Panov-def-extended} on the couple $(u^-,u^+)(t_0,\mx_0)$ of
one-sided traces of a solution $u$ at a point $(t_0,\mx_0)\in\Sigma$
where $\Sigma$ (called interface in the sequel) is a jump discontinuity hypersurface of a piecewise
smooth flux function $(t,\mx)\mapsto \mff(t,\mx,k)$ and $u^-,u^+$ are
strong one-sided traces of the admissible solution $u$ in the sense
of Definition \ref{def-adm}. The traces can be seen as limits, in
an appropriate sense, of $u((t_0,\mx_0)\pm h\nu|_\Sigma(t_0,\mx_0))$, as $h\downarrow 0$.
Here, $\nu|_\Sigma(t_0,\mx_0)$ is a normal vector to $\Sigma$ with some
orientation fixed by the choice of local coordinates in a neighbourhood of $
(t_0,\mx_0)$, and one can also define the traces of flux functions
 \begin{equation}\label{eq:1sidedfluxes}
 \mff^\pm(\cdot,k):=\lim_{h\downarrow 0}  \mff(\cdot\pm h\nu(\cdot),k).
 \end{equation}
The essence of the uniqueness argument is to justify the fact that in the
global Kato inequality  for two solutions $u$ and $v$, which formally reads
\begin{equation*}
\begin{array}{l}
  \pa_t |u-v|+\text{div}_\mx \Bigl({\rm sgn}(u-v)(\mff(t,\mx,u)-\mff(t,\mx,v))\Bigr)
\\
\qquad\qquad\qquad\qquad
\leq  \; L|u-v|\;+\; I\,  \mathcal H^{d}|_\Sigma \;\;{\rm in} \ \ {\cal D}'(\R^+\times\R^d),
\\[3pt]
I:=\Big\{{\rm sgn}(u^+\!-
v^+)(\mff^+(\cdot,u^+)-\mff^+(\cdot,v^+))\\
\qquad\qquad\qquad\qquad  - {\rm sgn}(u^-\!
-v^-)(\mff^-(\cdot,u^-)-\mff^-(\cdot,v_-))
\Big\}\cdot\nu(\cdot),
\end{array}
\end{equation*}
the contribution of the interface term $I$ on the right-hand side is non-positive.
Here $\mathcal H^{d}|_\Sigma$ is the $d$-dimensional Hausdorff measure supported on the interface $\Sigma$, and $L$ is the uniform in $(t,\mx)$ Lipschitz constant of the source term $S(t,\mx,\cdot)$. The non-positivity of $I$ follows indeed from the restrictions on $u^\pm$ imposed by \eqref{eq:Panov-def-extended}, see \S~\ref{ssec:FLAT-direct-uniq} (cf.~\cite{AKR}).

\smallskip
The genuine nonlinearity assumption on the flux, i.e., the assumption that
for every $(\xi_0,\xi)$ in the $d$-dimensional unit sphere $S^{d}$, there holds
\begin{equation}\label{gnl}
\begin{array}{l}
\text{for almost every $(t,\mx)\in \R^+\times \R^d$, the mapping}\\
k \mapsto \xi_0 k + \xi\cdot \mff(t,\mx,k)
\;\;\text{is not constant on non-degenerate intervals},
\end{array}
\end{equation}
is a technical hypothesis. It ensures existence of strong interface traces $u^\pm$ on $\Sigma$,
on which our uniqueness proof heavily relies. While this property can be circumvented in the space-time homogeneous
setting of \cite{AKR} and of many related references, in our case it is essential. Indeed,
our approach to existence of admissible solutions also relies on existence of $u^\pm$, but also on invariance of
the considered class of equations under variables' changes. Therefore, we have to work with $(t,\mx)$-dependent fluxes (and with source terms); in this general setting, assumption~\eqref{gnl}
is the only known condition that guarantees existence of traces (see \cite{AM_jhde} for more information on traces
of entropy solutions to an inhomogeneous conservation law).

\smallskip
Assumption \eqref{gnl} is essential for our existence proof also because it ensures compactness of suitable vanishing viscosity approximations.
It is standard, in this context, to assume in addition some conditions that ensure existence of invariant
regions for \eqref{eq:cons-law}; they are needed to obtain uniform $L^\infty$ bounds on sequences of
approximate solutions (see, e.g., \cite[Sect.\,6]{AKR}). Here, we take the simplest of such assumptions
that is:
\begin{equation}\label{eq:confinement}
 \begin{array}{l}
  \exists [a,b]\subset\R \;\;\text{such that for a.e. $t$, \; $\mff(t,\cdot,a)\equiv const$, $\mff(t,\cdot,b)\equiv const$,}
   \\
   \text{for a.e. $(t,\mx)$, one has $S(t,x,a)\geq 0$ and $S(t,x,b)\leq 0$,}\\
   \text{and the initial datum fulfills $a\leq u_0(\cdot)\leq b$}.
  \end{array}
\end{equation}
Assumption \eqref{eq:confinement} ensures that $u(t,\mx)\equiv a$ (respectively, $u(t,x)\equiv b$) is a sub-solution (resp., a super-solution) to the Cauchy problem for \eqref{eq:cons-law} with initial condition $u_0$, thus \eqref{eq:confinement} confines to the interval $[a,b]$ the values of solutions (and also the values of suitably constructed approximate solutions) we will consider.

\smallskip
We are able to prove existence of solutions defined by~\eqref{eq:Panov-def-extended},\eqref{eq:p-function}
  not by surpassing the technical difficulty of passage to the limit in approximate entropy
inequalities \eqref{eq:entropy-eps}  but using an indirect and rather lengthy argument.
 This argument exploits the idea of the existence proof in \cite{AKRnhm},
 developed in the simple case of a flat interface; see Remark~\ref{rem:Reps-1Dchoice} for an explanation of this choice. The study of the flat case is combined with appropriately chosen changes of variables and a principle of invariance of \eqref{eq:Panov-def-extended},\eqref{eq:p-function} under such variables' changes (see in particular Proposition~\ref{prop:equiv-changeofvar}).
We are able to include changes of variables that are singular
in a neighbourhood of $(d-1)$-dimensional singularities such as intersections of different $d$-dimensional interfaces.
As a typical and illuminating example, we consider the case of interface $\Sigma\subset \R^+\times\R$
consisting in two Lipschitz curves merging into one.


\subsection{Outline of the paper}\label{ssec:INTRO-outline}
The paper is organized as follows.
In Section~\ref{sec:FLAT} we treat the
case that is fundamental for our techniques.  More precisely,
in this section we deal with multi-dimensional scalar
conservation laws in heterogeneous setting, with a source term but with a single flat
discontinuity interface $\Sigma$. Definition of solution admissibility is stated in \S~\ref{ssec:FLAT-def+traces}
and further discussed in \S~\ref{ssec:FLAT-pu-specified}. Uniqueness of an admissible solution to the Cauchy problem is proved in \S~\ref{ssec:FLAT-direct-uniq}, then a different proof expoiting the machinery of \cite{AKR,AKRnhm} is sketched in \S~\ref{ssec:FLAT-GVV-uniq}. Existence is justified in \S~\ref{ssec:FLAT-exist}, with an explicit construction of
$p_u$ in terms of interface traces $u^\pm$ of the vanishing viscosity limit $u$.

\smallskip
Section~\ref{sec:CURVED} deals with the general geometry of jump discontinuity interfaces of the flux $\mff$ in \eqref{eq:cons-law}. In \S~\ref{ssec:CURVED-prelim}, we sketch the uniqueness argument, explain and motivate the existence strategy. In \S~\ref{ssec:CURVED-AlmostRect}, we give main details on the idea of a global singular
change of variables that permits to reduce the case of flux $\mff $ with Lipschitz jump discontinuity manifolds
of rather general form (possibly curved, intersecting, and closed) to the locally flat case with singularities that can be ignored. With minor modifications with respect to \S~\ref{ssec:FLAT-exist}, existence in this locally flat case is justified. In \S~\ref{ssec:CURVED-rectifiable-examples}, examples of changes of variables satisfying the general assumptions of \S~\ref{ssec:CURVED-AlmostRect} are presented. In addition, in \S~\ref{ssec:CURVED-LocalDiffusion} we briefly describe an alternative
construction of solutions where regular local changes of variables are pieced together using a kind of partition of unity.
In Remark~\ref{rem:viscosity-bis}, we underline the fact that in both approaches, a viscosity approximation
adapted to interface geometry is essential for the existence proof.

\smallskip
For readers' convenience, let us point out that Definition~\ref{def-adm} with well-posedness results of Theorems~\ref{th:uniq},\ref{thm5} summarize our work for the model case with flat interface.
For the general case of piecewise Lipschitz, jump discontinuous in $(t,\mx)$ flux $\mff$, Theorem~\ref{th:well-posedness-general} (see also the conclusion of \S~\ref{ssec:CURVED-LocalDiffusion}) provides
 a set of sufficient conditions under which the new interpretation \eqref{eq:Panov-def-extended},\eqref{eq:p-function} of equation~\eqref{eq:cons-law} leads to well-posedness of the Cauchy problem.

\section[New admissibility conditions for flat discontinuity]{New admissibility conditions for multidimensional heterogeneous scalar conservation laws with a flat discontinuity}\label{sec:FLAT}

This is the fundamental section, indeed, all the other situations that we can
resolve are reduced, at least locally, to the Cauchy problem investigated here.

\subsection{Definition of solution and interface traces}\label{ssec:FLAT-def+traces}
Consider the problem
\begin{align}
\label{flat} &\pa_t u + \Div (F(t,\mx,u)H(x_1)+F(t,\mx,u)H(-x_1))=S(t,\mx,u),\\
\label{flati}&u|_{t=0}=u_0, 
\end{align}
where $F=(F_1,\dots,F_d):\R^+\times\R^d\times [a,b]\to \R^d$ and $G=(G_1,\dots,G_d):\R^+\times\R^d\times [a,b]\to
\R^d$.
This is \eqref{eq:cons-law} with the flux given by
\begin{equation}\label{eq:flat-flux}
\mff(t,\mx,\cdot)=F(t,\mx,\cdot)H(x_1)+G(t,\mx,\cdot)H(-x_1).
\end{equation}
 We will assume that the flux and the source $S$ satisfy the compatibility conditions at $k=a$ and $k=b$,
as required in assumption \eqref{eq:confinement}. We also assume that
\begin{equation}\label{eq:Lip-flux}
\begin{array}{l}
  \text{the fluxes $F,G$ are globally Lipschitz continuous in all variables,}\\
  \text{the source $S$ is globally Lipschitz continuous in $u\in [a,b]$}.
\end{array}
\end{equation}
Here and in the sequel of this paper,
$\hat{x}_1=(x_2,\dots,x_d)$. The unique interface (jump discontinuity hypersurface) for the flux~\eqref{eq:flat-flux}
is given by
$$\Sigma=\{(t_0,0,\hat{x}_{10})\,|\,t_0\in\R^+,\hat x_{10}\in\R^{d-1}\}.$$
Throughout this section, somewhat abusively we denote by $(t_0,\mx_0)$ both the points of $\Sigma$ with coordinates $(t_0,\hat x_{10})$ and the points of $\R\times\R^{d}$ with coordinates $(t_0,0,\hat x_{10})$.

\smallskip
As explained in the introduction, the main goal of this paper is to propose and justify, by proving well-posedness results,
the following definition.
\begin{definition}
\label{def-adm} We say that a function $u\in
C(\R^+;L^1_{loc}(\R^d))$ taking values in $[a,b]$ is an entropy
admissible solution to \eqref{flat}, \eqref{flati} if
$u(0,\cdot)=u_0$ and
 there exists a function $p_u:[0,T]\to [a,b]$ such that $\forall k\in
 [a,b]$
\begin{align}
\label{sep827} & \pa_t |u\!-\!k|+\Div_\mx
\Big\{{\rm sgn}(u\!-\!k)\Bigl(H(x_1)({F}(t,\mx,u)\!-\!{F}(t,\mx,k))\\
&\qquad\qquad+H(-x_1)({G}(t,\mx,u)\!-\!{G}(t,\mx,k))\Bigr)\Big\} \;-\;\text{sgn}(u-k)S(t,\mx,u) \nonumber\\
&+{\rm sgn}(u-k) \Bigl(H(x_1)\Div_\mx
F(t,\mx,k)+H(-x_1)\Div_\mx G(t,\mx,k)\Bigr) \nonumber\\
&  -{\rm
sgn}(p_u(t,\hat{x}_1)\!-\!k)({F_1}(t,\mx,k)\!-\!{G_1}(t,\mx,k))\delta_0(x_1) \;\leq\; 0
 \ \ {\rm in} \ \ {\cal D}'(\R^+\times \R^d)\nonumber.
\end{align}
\end{definition}
To shorten the calculations, in the sequel we will assume $S\equiv 0$; the general case is
obtained with the help of the Gronwall inequality. The technique of the proofs of this section
readily extends to a locally finite number of flat, possibly crossing discontinuity surfaces, but we will stick to the case of one interface $\Sigma=\{x_1=0\}$.
Let us stress that presence of a source term $S$ and of multiple flat
discontinuity hypersurfaces is required in order to reduce the case of curved, crossing or closed interfaces considered in Section~\ref{sec:CURVED} to the case of flat interfaces investigated in the present section.

Clearly,  inequalities \eqref{sep827}  (with $S\equiv 0$) in the one-dimensional
situation imply the Karlsen-Risebro-Towers inequalities \eqref{eq:KRT-cond}, therefore it is natural
to expect that the uniqueness
becomes easier to justify, while existence proof will present new and considerable difficulties.
 Indeed, in the present section we will show that
\begin{itemize}
  \item an admissible solution in the sense of \eqref{sep827} is unique,
  with a rather standard proof involving a tedious case-by-case study;

  \item the standard vanishing viscosity method converges towards this solution,
  with a quite indirect proof based upon construction of viscosity profiles.
\end{itemize}
Both results are achieved by looking at the values of one-sided traces of the solutions on the interface,
so we start by defining these traces and giving sufficient conditions for their existence.

%

\begin{definition}
 We say that an integrable function $W$  admits left trace $W^-$ and right trace $W^+$
 on $\Sigma=\{x_1=0\}$ if there exist functions $W^\pm: \R^+\times\R^{d-1}\mapsto \R$ such that
 for every $\varphi\in C_c(\R^{d-1})$ there holds
\begin{equation}\label{eq:traces}
 \lim\limits_{\substack{h\to 0^+}} \int_{\R^+\times\R^{d-1}}
|W(t,\pm h,\hat{x}_1)-W^\pm(t,\hat{x}_1)| \varphi(t,\hat{x}_1)
dt d \hat{x}_1 = 0.
\end{equation}
\end{definition} While general $L^\infty$ functions do not admit such traces, local
entropy solutions of conservation laws admit traces in the sense
\eqref{eq:traces} under additional, not very restrictive technical
assumptions. Under assumption \eqref{gnl} on the flux~\eqref{eq:flat-flux},
 the
arguments from \cite{pan_jhde} provide existence of traces to
local Kruzhkov entropy solutions to \eqref{flat} in the non-homogeneous situation as
well (see also \cite{AM_jhde}). Indeed, while the proof from \cite{pan_jhde} is given
for case where the flux does not depend on time or space variables
(i.e. it is homogeneous), its arguments extend to the general,
piecewise Lipschitz $(t,\mx)$-dependent setting with condition
\eqref{gnl} (see in particular \cite{AM_jhde}).
 Let us mention in passing that in one
dimensional homogeneous situation, the linear degeneracy of fluxes
can also be treated via introduction of ``singular mappings'' that
have traces in all cases (see \cite{AKR}) or by considering traces
of the flux and of entropy fluxes instead of the traces of the
solution itself (see \cite{dm_dcds}).

Throughout this section, we will denote by $u^-(t_0,\mx_0)$, respectively $u^+(t_0,\mx_0)$, the left trace, respectively the right trace at $(t_0,\mx_0)\in\Sigma$ of an admissible solution $u$ to \eqref{flat}.
 We readily derive the Rankine-Hugoniot conditions for admissible solutions.
\begin{lemma}\label{lem:RH}
Let $u$ be an admissible solution to \eqref{flat}, \eqref{flati} in the which
admits left and right strong traces at the interface $u^+$ and
$u^-$, respectively. Then, for $\mathcal H^d$-a.e. point $(t_0,\mx_0)\in \Sigma$, there
holds
\begin{equation}
\label{R-H}
 F_1(t_0,\mx_0,u^+(t_0,\mx_0))= G_1(t_0,\mx_0,u^+(t_0,\mx_0))
\end{equation}
\end{lemma}
\begin{proof}
Taking $k =a$ then $k=b$ in the entropy inequalities \eqref{sep827}, one deduces the weak formulation of \eqref{flat}.
Then it is enough to test this weak formulation of \eqref{flat} on the function of the form $\eta=
\varphi \mu_h$, $\varphi\in C^1_c(\R^+\times \R^d)$, where $\mu_h$ is
given by
\begin{equation}
\label{mih} \mu_h(\mx)=
\begin{cases}
\frac{1}{h}(x_1+2h), & x_1\in [-2h,-h]\\
1, & x\in [-h,h]\\
\frac{1}{h}(2h-x_1), & x_1\in [h,2h]\\
0, & |x_1|>2h,
\end{cases}
\end{equation}
which can be taken for test function by approximation.
 After letting $h\to 0$, we reach to identity \eqref{R-H} in $\mathcal D'(\Sigma)$; since both sides of \eqref{R-H}
 are bounded functions, the identity also holds pointwise, $\mathcal H^d$-a.e. on $\Sigma$.
\end{proof}

\subsection{A uniqueness and $L^1$ contraction proof}\label{ssec:FLAT-direct-uniq}
We prove the following result. Let $B(0,R)=\{x\in\R^d\,|\, |x|<R\}$.
\begin{theorem}\label{th:uniq}
Assume that the flux \eqref{eq:flat-flux} satisfies
\eqref{gnl},\eqref{eq:confinement},\eqref{eq:Lip-flux}. Let $u,v$
be two admissible solutions to \eqref{flat} with source $S\equiv 0$ and initial data $u_0$ and $v_0$, respectively.
Then, for every
$R>0$, $T>0$ there exists $C>0$ such that
\begin{equation}
\label{stab} \int_0^T\int_{B(0,R)}|u(t,x)-v(t,x)|dxdt \leq T
\int_{B(0,R+CT)}|u_0(x)-v_0(x)|dx.
\end{equation}
\end{theorem}
Notice that the order-preservation result ($u_0\geq v_0$ implies $u\geq v$) can be proved in the same way.
The first proof given below is a self-contained case-by-case study.
The second proof that will be sketched in \S\,\ref{ssec:FLAT-GVV-uniq} establishes the equivalence between
our new definition and the definition of $\mathcal G_{VV}$-entropy solution (see \cite{AKR,AKRnhm}) for which
uniqueness follows readily from the general theory developed in \cite{AKR} (see also \cite{Die3}).

\smallskip
We start as in the proof of Lemma~\ref{lem:RH}, but using the entropy inequalities with general $k$. Namely, insert into \eqref{sep827} the function $\psi=\mu_h \varphi$, where $\varphi \in
C^1_c(\R^d_+)$ while $\mu_h$ is given by \eqref{mih}. Letting $h\to 0$, due to
arbitrariness of $\varphi$, as in Lemma~\ref{lem:RH} we reach to the following relation for
almost every $(t_0,\mx_0)\in \Sigma$
\begin{equation}
\label{sw-cond}
\begin{split}
&{\rm sgn}(u^+-k)(F_1(t_0,\mx_0,u^+)-F_1(t_0,\mx_0,k))\\&
\quad -{\rm
sgn}(u^--k)(G_1(t_0,\mx_0,u^-)-G_1(t_0,\mx_0,k))\\&
\qquad +{\rm
sgn}(p_u-k)(F_1(t_0,\mx_0,k)-G_1(t_0,\mx_0,k))
\leq 0,
\end{split}
\end{equation}
where for the sake of brevity, we write
 $$
p_u(t_0,\mx_0)=p_u, \ \ u^+(t_0,\mx_0)=u^+, \ \
u^-(t_0,\mx_0)=u^-.
$$
Notice that in the passage to the limit, we used continuity of the maps
\begin{equation}\label{eq:continuity-of-KrFluxes}
 \begin{array}{l}
   (t,\mx,z)\mapsto \text{sgn}(z-k)(F(t,\mx,z)-F(t,\mx,k)\\
 (t,\mx,z)\mapsto \text{sgn}(z-k)(G(t,\mx,z)-G(t,\mx,k).
 \end{array}
\end{equation}

 Remark that if
$$ \text{  $k \geq \max\{u^+, u^-,
p \}$ \;or\; $k \leq \min\{u^+,
u^-, p \}$
}
$$
 then the left-hand side in \eqref{sw-cond} is equal to zero according to the Rankine-Hugoniot relation \eqref{R-H}. Now, we make precise the information contained in the other possible cases where \eqref{sw-cond} holds true.
To shorten the statements like~\eqref{R-H},\eqref{sw-cond} in the sequel we will use the notation
 \begin{equation}\label{eq:f-and-g}
   \text{$f(\cdot):=F_1(t_0,\mx_0,\cdot)$ and $g(\cdot):=G_1(t_0,\mx_0,\cdot)$}.
 \end{equation}

\medskip
\noindent {\bf Case I:} $ \ \ u^+\leq u^-$.

\smallskip
\noindent$\bullet$
 $u^+\leq u^-\leq p_u$, where two different cases occur:
\begin{equation}
\label{i1}
u^+\leq u^-\leq k \leq p_u \, , \;\;\text{which implies}\;\;
f(k)\leq g(k);
\end{equation}
%
\begin{equation}
\label{ii1}
u^+\leq k \leq u^- \leq p_u \, , \;\;\text{which implies}\;\;
f(k)\leq f(u^+).
\end{equation}

\smallskip
\noindent$\bullet$ $u^+\leq p_u \leq u^-$, where we also have two possible cases:
\begin{equation}
\label{iii1}
u^+ \leq p_u \leq k \leq u^- \, , \;\;\text{which implies}\;\;
g(k)\leq g(u^-);
\end{equation}
%
\begin{equation}
\label{iv1}
u^+\leq k \leq p_u \leq u^-  \, , \;\;\text{which implies}\;\;
f(k)\leq f(u^+).
\end{equation}

\smallskip
\noindent$\bullet$ $p \leq u^+\leq u^-$, here we have the following alternative:
\begin{equation}
\label{v1}
p \leq u^+ \leq  k \leq u^- \, , \;\;\text{which implies}\;\;
f(k)\leq f(u^+);
\end{equation}
\begin{equation}
\label{vi1}
p \leq  k \leq u^+ \leq   u^-  \, , \;\;\text{which implies}\;\;
g(k)\leq f(k).
\end{equation}

\medskip
\noindent {\bf Case II:} $ \ \ u^-\leq u^+$.

\smallskip
\noindent$\bullet$ $u^-\leq u^+\leq p_u$, where we have one of the two situations:
\begin{equation}
\label{i2}
u^-\leq u^+\leq k \leq p_u \, , \;\;\text{which implies}\;\;
f(k)\leq g(k);
\end{equation}
\begin{equation}
\label{ii2}
u^-\leq k \leq u^+ \leq p_u \, , \;\;\text{which implies}\;\;
g(u^-)\leq g(k).
\end{equation}

\smallskip
\noindent$\bullet$
 $u^-\leq p_u \leq u^+$, where the alternative is:
\begin{equation}
\label{iii2}
u^- \leq p_u \leq k \leq u^+ \, , \;\;\text{which implies}\;\;
f(u^+)\leq f(k);
\end{equation}
%
\begin{equation}
\label{iv2}
u^-\leq k \leq p_u \leq u^+  \, , \;\;\text{which implies}\;\;
g(u^-)\leq g(k).
\end{equation}

\smallskip
\noindent$\bullet$
 $p_u \leq u^-\leq u^+$, where we have the last two possibilities:
\begin{equation}
\label{v2}
p_u \leq u^- \leq  k \leq u^+ \, , \;\;\text{which implies}\;\;
f(u^+)\leq f(k);
\end{equation}
%
\begin{equation}
\label{vi2}
p_u \leq  k \leq u^- \leq   u^+  \, , \;\;\text{which implies}\;\;
g(k)\leq f(k).
\end{equation}

\smallskip
Now, we are ready to give a proof of the uniqueness result of Theorem~\ref{th:uniq}.
\begin{proof}
The main part of the proof consists in derivation of the Kato inequality:
for every $\varphi\in C^1_c(\R^+\times\R^d)$,
\begin{align}
\label{kato} \int_{\R^+\times\R^d}\Bigl\{& |u-v|\varphi_t+{\rm
sgn}(u-v)\Bigl(
(F(t,\mx,u)-F(t,\mx,v))H(x_1)\\
&+(G(t,\mx,u)-G(t,\mx,v))H(-x_1)\Bigr) \cdot\nabla \varphi \Bigr\}d\mx dt
\geq 0. \nonumber
\end{align}
The classical doubling of variables technique of \cite{Kru} ensures that \eqref{kato} holds for $\varphi\in
C^1_c((\R^+\times\R^d)\setminus \Sigma)$.
Therefore, given an arbitrary $\psi\in C^1_c(\R^+\times \R^d)$, inequality \eqref{kato}
is satisfied with the test function $\varphi= (1-\mu_h)\psi$, where $\mu_h$ is defined by \eqref{mih}. Letting
$h\to 0$, we get
\begin{align}
\label{kato-per} &\int_{\R^+\times\R^d}\!\!\Bigl\{ |u\!-\! v|\psi_t
\!+ \! {\rm sgn}(u\!-\! v)\Bigl(
(F(t,\mx,u)\!-\!F(t,\mx,v))H(x_1)\!\\&\qquad\qquad\qquad+\!(G(t,\mx,u)\!-\!G(t,\mx,v))H(-x_1)\Bigr)
\nabla \psi \Bigr\}d\mx dt \nonumber\\
&\geq \int_{\Sigma}\!\Bigl( -{\rm
sgn}(u^+\!-\!v^+)(F(t_0,\mx_0,u^+)\!-\!F(t_0,\mx_0,v^+))\!\nonumber\\&\qquad\qquad\qquad+\!{\rm
sgn}(u^-\!-\!v^-)(G(t_0,\mx_0,u^-)\!-\!G(t_0,\mx_0,v^-))\Bigr)\psi(t_0,\mx_0)\,d\hat{x}_{10}
dt_0\nonumber\\&=:
\int_{\Sigma} \Delta(u^\pm,v^\pm)(t_0,\mx_0)\psi(t_0,\mx_0) d\hat{x}_{10}\, dt_0,\;\; \mx_0=(0,\hat x_{10}).
\nonumber
\end{align}
 Now, we prove
that the integrand  $\Delta(u^\pm,v^\pm)$ (that is a short-cut notation for $\Delta(t_0,\mx_0,u^\pm(t_0,\mx_0),v^\pm(t_0,\mx_0))$) in the right-hand side of the latter expression is greater than or equal to zero, for almost every $(t_0,\mx_0)=(t,0,\hat{x}_{10})\in \Sigma$. The proof is tedious and it is accomplished by
considering numerous different possibilities depending on pointwise relations
between the values $u^\pm$, $v^\pm$ and $p_u$, $p_v$.

Concerning the relation between $u^\pm$ and $v^\pm$, we see that,
according to the Rankine-Hugoniot conditions, the quantity $\Delta (u^\pm,v^\pm)$ can be non-zero
 only when
\begin{equation}
\label{crit}
u^->v^- \  {\rm and}  \  u^+<v^+  \ \ \ {\rm or} \ \ \
u^-<v^- \  {\rm and}  \ u^+>v^+,
\end{equation}
  (the proof is the same as in Cases 1--5 in the proof of \cite[Theorem 2.1]{kar3}).
  On the other hand, the two cases from \eqref{crit} are symmetric since
$\Delta(u^\pm,v^\pm)=\Delta(v^\pm,u^\pm)$, and therefore their analysis is
the same. Thus, it is enough to prove the inequality $\Delta(u^\pm,v^\pm)
\geq 0$ whenever the first relation from \eqref{crit} is satisfied.
We proceed by considering the following possible sub-cases:
\begin{align*}
&\text{Case 1:} \qquad u^+ < v^+ <v^-<u^- \\
&\text{Case 2:} \qquad u^+ < v^- <v^+<u^- \\
&\text{Case 3:} \qquad v^- < u^+ <v^+<u^- \\
&\text{Case 4:} \qquad v^- < u^+<u^-<v^+\\
&\text{Case 5:} \qquad v^- < u^-<u^+<v^+.
\end{align*}
Notice that, according to the disposition of $u^\pm$
and $v^\pm$ and the Rankine-Hugoniot condition~\eqref{R-H}, one has
\begin{equation*}
\begin{split}
&\Delta(u^\pm,v^\pm)=-{\rm sgn}(u^+\!-\!v^+)(f(u^+)\!-\!f(v^+)) +
{\rm sgn}(u^-\!-\!v^-)(g(u^-)\!-\!g(v^-))\\
&\quad =
f(u^+)-f(v^+)+g(u^-)-g(v^-)=2(f(u^+)-f(v^+))=2(g(u^-)-g(v^-)).
\end{split}
\end{equation*} Thus, we aim to prove that in each case of the above list, there holds
\begin{equation}
\label{aim} f(u^+)-f(v^+)\geq 0 \ \ \
{\rm or} \ \ \ g(u^-)-g(v^-) \geq 0.
\end{equation}

Since $u,v\in L^\b(\R^+\times\R^d)$ are two admissible solutions to
\eqref{flat}, we consider two functions $p_u$ and
$p_v$ from Definition \ref{def-adm} representing the singular values on $\Sigma$
corresponding to $u$ and $v$, respectively.

\medskip
\noindent
{\bf Case 1} $ \ \ $ For almost every fixed $(t_0,\mx_0) \in \Sigma$,
we have the following possibilities.

\smallskip
\noindent$\bullet$ $ \ \ u^+ < v^+ <v^-<u^-<p_u$:

The conclusion follows by taking $k=v^+$ in \eqref{ii1}.

\noindent$\bullet$ $ \ \ u^+ < v^+ \leq p_u \leq u^-$:

The conclusion follows by taking $k=v^+$ in \eqref{iv1}.

\noindent$\bullet$ $ \ \ u^+ \leq p_u \leq v^- <u^-$:

The conclusion follows by taking $k=v^-$ in \eqref{iii1}.

\noindent$\bullet$ $ \ \ p_u< u^+ < v^+ <v^-<u^-$:

The conclusion follows by taking $k=v^-$ in \eqref{v1}.

\medskip\noindent
{\bf Case 5} $ \ \ $ This case is symmetric with the previous one.
We simply need to  consider the position of $p_v$ instead of $p_u$ and to
apply \eqref{i2}--\eqref{vi2} instead of \eqref{i1}--\eqref{vi1}.

\medskip\noindent
{\bf Case 2} $ \ \ $ We have the
following possibilities.

\smallskip
\noindent$\bullet$ $ \ \ u^+ < v^- <v^+<u^-<p_u$:

The conclusion follows by taking $k=v^+$ in \eqref{ii1}.

\noindent$\bullet$ $ \ \ u^+ < v^-<v^+ \leq p_u \leq u^-$:

The conclusion follows by taking $k=v^+$ in \eqref{iv1}.

\noindent$\bullet$ $ \ \ u^+ < v^-\leq p_u \leq v^+ < u^-$:

Here, we must involve the position of $p_v$. First, recall
that from \eqref{iii1} and \eqref{iv1}
\begin{equation}
\label{help1}
\begin{split}
&g(k)\leq g(u^-), \ \ \forall k \in [p_u,u^-]\\
&f(k)\leq f(u^+), \ \ \forall k \in [u^+,p_u].
\end{split}
\end{equation}Now, we have the
following possibilities.
\begin{itemize}
\item[$\cdot$]  $ \ \ u^+ < v^-\leq p_u \leq v^+ \leq p_v$.\\
From \eqref{ii2} (applied on $v$) and \eqref{help1}, we have
 $g(v^-)\leq g(p_u) \leq g(u^-)$ which is \eqref{aim}.

\item[$\cdot$]  $ \ \ u^+ < v^-\leq p_u \leq p_v \leq v^+<u^- $.\\
From \eqref{iv2} and \eqref{help1}, we have $g(v^-)\leq
g(p_u)\leq g(u^-)$.

\item[$\cdot$]  $ \ \ u^+ < v^-\leq p_v \leq p_u \leq v^+<u^- $.\\
From \eqref{iii2} and \eqref{help1}, we have
$f(v^+)\leq f(p_u)\leq f(u^+)$.

\item[$\cdot$]  $ \ \ p_v \leq v^-\leq p_u \leq  v^+<u^- $.\\
From \eqref{v2} and \eqref{help1}, we have $f(v^+)\leq
f(p_u)\leq f(u^+)$.
\end{itemize}

\noindent$\bullet$ $ \ \ u^+ \leq p_u \leq v^-< v^+ <u^-$:

The conclusion follows by taking $k=v^-$ in \eqref{iii1}.

\noindent$\bullet$ $ \ \ p_u< u^+ < v^+ <v^-<u^-$:

The conclusion follows by taking $k=v^-$ in \eqref{v1}.

\medskip\noindent
{\bf Case 4} $ \ \ $ This case is symmetric with the previous one.
We simply need to consider the position of $p_v$ instead of $p_u$ and to
apply \eqref{i2}--\eqref{vi2} instead of \eqref{i1}--\eqref{vi1} or
vice versa, when needed.

\medskip\noindent
{\bf Case 3} $ \ \ $ We have the
following possibilities.

\smallskip
\noindent$\bullet$ $ \ \ v^-<u^+<v^+<u^-\leq p_u$:

In this case, the first relation in \eqref{aim} follows by taking
$k=v^+$ in \eqref{ii1}.

\noindent$\bullet$ $ \ \ v^-<u^+<v^+\leq p_u\leq u^-$:

In this case, \eqref{aim} follows from \eqref{iv1} by taking
$k=v^+$ there.

\noindent$\bullet$ $ \ \ v^-<u^+\leq p_u \leq v^+< u^-$:

We must involve the position of $p_v$ again. We have the following
possibilities.
\begin{itemize}
\item[$\cdot$] $\ \ v^-<u^+\leq p_u \leq v^+< u^-\leq p_v$.\\
From \eqref{ii2}, on $v$, it follows $g(v^-)\leq g(p_u)$ while from
\eqref{help1}, $g(p_u)\leq g(u^-)$. Thus, \eqref{aim} follows.

\item[$\cdot$]  $\ \ v^-<u^+\leq p_u \leq v^+\leq p_v \leq  u^-$.\\
The situation is the same as the previous one.

\item[$\cdot$]  $\ \ v^-<u^+\leq p_u \leq p_v \leq v^+ \leq  u^-$.\\
From \eqref{iv2}, on $v$, and \eqref{help1}, it follows
$g(v^-)\leq g(p_u) \leq g(u^-)$ which is \eqref{aim}.

\item[$\cdot$]  $\ \ v^- < u^+ \leq p_v \leq p_u \leq v^+ \leq  u^-$.\\
From \eqref{iii2} and \eqref{help1}, it follows
$f(v^+)\leq f(p_u) \leq f(u^+)$ which is \eqref{aim}.

\item[$\cdot$]  $\ \ v^- \leq p_v \leq  u^+  \leq p_u \leq v^+ \leq  u^-$.\\
Relation \eqref{aim} follows as in the previous case.

\item[$\cdot$]  $\ \ p_v \leq v^- <   u^+  \leq p_u \leq v^+ \leq  u^-$.\\
From \eqref{v2} and \eqref{help1}, it follows
$f(v^+)\leq f(p_u) \leq f(u^+)$.
\end{itemize}

\noindent$\bullet$  $ \ \ v^-\leq p_u \leq u^+ < v^+< u^-$:

Relation \eqref{aim} follows from \eqref{v1}.

\noindent$\bullet$  $ \ \ p_u \leq v^-<  u^+ < v^+< u^-$:

The conclusion is the same as in the previous item.

\smallskip
From the above considerations, we conclude that in all possible cases, there holds $\Delta(u^\pm,v^\pm)\geq
0$. This means that the Kato inequality~\eqref{kato} holds.
From here, the proof of the theorem follows in the standard way, see \cite{Kru,AKR}.
\end{proof}

\subsection{Uniqueness via reduction to $\mathcal G_{VV}$-entropy solutions}\label{ssec:FLAT-GVV-uniq}

It is possible to reduce uniqueness proof to the setting of \cite{AKR,AKRnhm}.
Indeed, in these papers uniqueness of so-called $\mathcal G_{VV}$-entropy solutions has been proved.
Here, given an ordered couple of continuous functions $(f,g)$ on $[a,b]$, the \emph{vanishing viscosity germ} $\mathcal G_{VV}$ is the subset of $[a,b]^2$ given by: $(u^-,u^+)\in \mathcal G_{VV}$ if $g(u^-)=f(u^+)=:s$ and
  \begin{equation}\label{eq:ClosureViscGerm}
    \begin{array}{l}
      \hspace*{20pt}\text{\emph{either} $u^-=u^+$,} \\
      \hspace*{20pt}\text{\emph{or} $u^-<u^+$ and there exists}\\
      \hspace*{20pt}\hskip 1cm \text{$u^o\in \left[u^-,u^+\right]$ such that}
      \begin{cases}
        g(z)\ge s\ \text{for all $z\in \left[u^-,u^o\right]$,}\\
       f(z)\ge s\ \text{for all $z\in\left[u^o,u^+\right]$,}
      \end{cases}\\
      \hspace*{20pt}\text{\emph{or} $u^->u^+$ and there exists}\\
      \hspace*{20pt}\hskip 1cm \text{$u^o\in \left[u^-,u^+\right]$ such that}
      \begin{cases}
        g(z)\le s\ \text{for all $z\in \left[u^o,u^-\right]$,}\\
        f(z)\le s\ \text{for all $z\in\left[u^+,u^o\right]$.}
      \end{cases}
    \end{array}
  \end{equation}
  This set is called the \emph{vanishing viscosity germ} associated
  with the couple $(f,g)$. We refer to Diehl in \cite{Die3} for an equivalent description of the set $\mathcal G_{VV}$ in terms of Oleinik-kind inequalities. Then, $\mathcal G_{VV}$-entropy solutions are defined
  as follows.
\begin{definition}
\label{gvv-adm}We say that the function $u\in L^\infty(\R^+\times
\R^d)$ represent $\mathcal G_{VV}$-entropy solution to \eqref{flat},
\eqref{flati} if

\noindent\smallskip
$\bullet$ \eqref{sep827} holds when tested on functions $\varphi\in
C^1_c((\R^+\times \R^d)\setminus \Sigma)$;

\noindent\smallskip
$\bullet$  for a.e. $(t_0,\mx_0)\in\Sigma$, left $u^-(t_0,\mx_0)$ and right $u^+(t_0,\mx_0)$ traces of $u$ on $\Sigma$  are such that $(u^-(t_0,\mx_0),u^+(t_0,\mx_0))$ belongs to the vanishing viscosity germ associated with the couple
$(f,g)$ given by $f:u\mapsto F_1(t_0,\mx_0,u)$, $g:u\mapsto G_1(t_0,\mx_0,u)$.
\end{definition}
Then it is not difficult to check that a solution is admissible in the sense of Definition~\ref{def-adm}
if and only if it is a $\mathcal G_{VV}$-entropy solution. Indeed, one can take for $p_u(t_0,\mx_0)$
in \eqref{sep827} the value $u^o(t_0,\mx_0)$ associated with the couple $(u^-,u^+)(t_0,\mx_0)$ in definition~\eqref{eq:ClosureViscGerm} of $\mathcal G_{VV}$, and vice versa. The complete verification
 involves essentially
the same case studies as in the direct uniqueness proof developed in the previous section, therefore we omit these details.

\subsection{The existence proof}\label{ssec:FLAT-exist}

Recall that we assumed that the confinement property~\eqref{eq:confinement} holds, and that equation \eqref{flat} is genuinely nonlinear in the sense \eqref{gnl}. These properties imply that the family of vanishing viscosity approximations for \eqref{flat},\eqref{flati} (i.e., solutions of \eqref{3} below with initial condition~\eqref{flati})
is $[a,b]$-valued and that it is strongly precompact  (\cite{pan_arma}); moreover, the corresponding limit admits
strong traces at the interface $\Sigma$ (cf. \cite{AM_jhde}). Recall that we assume that $F,G$ are given Lipschitz functions on $\R^+\times\R^d\times[a,b]$. In addition, for the proof of existence we assume
that
\begin{equation}\label{hyp:pa_u F,G continuous in t,x}
\text{$\pa_u F$, $\pa_u G$ are Lipschitz continuous in $(t,x)$, $\mathcal H^d$-a.e. on $\Sigma$.}
\end{equation}
This rather strong, but not restrictive in practice assumption
simplifies the proof below; mere continuity and even a uniform in
$k$ Lebesgue-point property of $\pa_u F(\cdot,\cdot,k)$, $\pa_u
G(\cdot,\cdot,k)$ at $\mathcal H^d$-a.a. point of $\Sigma$ would be
enough for the technique of the proof to work.

Let $u$ be constructed as an accumulation point, as $\eps\to 0$, of vanishing viscosity approximations $(u_\eps)_\eps$. The essence of the existence proof consists in construction of a function $p_u$ on $\Sigma$ in order to justify that $u$ is an admissible solution of \eqref{flat},\eqref{flati} in the sense of Definition~\ref{def-adm}.
We will need several auxiliary statements before the conclusion can be given.
First, we have two lemmas that permit to extract information on $u=\lim_{\eps\to 0} u_\eps$ from
existence of a suitable family of vanishing viscosity profiles $(R_\eps)_\eps$ with prescribed traces of $R:=\lim_{\eps\to 0} R_\eps$ on $\Sigma$. We start by justifying a Kato inequality for $u_\eps$ and $R_\eps$, where $R_\eps$ solves an auxiliary equation with frozen coefficients. Recall that we assumed $S=0$.
\begin{lemma}
\label{l1} Assume that in a neighbourhood $V(t_0,\mx_0)$ of a point $(t_0,\mx_0)\in\Sigma$,
the family of functions $(u_\eps)_\eps$ in $L^2(0,T;H^1_{loc}(\R^d))$ satisfies
equations
\begin{equation}
\label{3} \pa_t
u_\eps+\Div_\mx \left(H(x_1)F(t,\mx,u_\eps)+H(x_1)G(t,\mx,u_\eps)
\right)=\eps \Delta u_\eps
\end{equation}
in the sense of distributions.
Assume that the family of  $L^2(0,T;H^1_{loc}(\R^d))$ functions $(R_\eps)_\eps$ takes values in $[a,b]$,
has its variation is uniformly bounded
by a constant $M$, and $R_\eps$ satisfies
\begin{equation}
\label{4}
\pa_t R_\eps \!+ \! \Div_\mx
\left(F(t_0,\mx_0,R_\eps)H(x_1)
+G(t_0,\mx_0,R_\eps)H(-x_1)\right)=\eps \Delta R_\eps
\end{equation} in the same neighbourhood $V(t_0,\mx_0)$ of $(t_0,\mx_0)$ in the sense of distributions.

Then, the following Kato-type inequality is satisfied:
\begin{align}
\label{5} &\pa_t|u_\eps -R_\eps| +\Div_x \Big( {\rm sgn}(u_\eps
\!-\!R_\eps) ((F(t,x,u_\eps)\!-\!F(t,x,R_\eps))H(x_1)\!
\\&\quad  +\!
(G(t,x,u_\eps)\!-\!G(t,x,R_\eps))H(-\!x_1)) \Big)\;  \leq \; \gamma_\eps(t,\mx)+\eps \Delta |u_\eps-R_\eps| \nonumber
\end{align}
in the sense of distributions in $V(t_0,\mx_0)$,
 where $(\gamma_\eps)_\eps$ 
  is a family of
 Radon measures in $V(t_0,\mx_0)$ verifying
\begin{equation}\label{eq:bound-on-gamma}
  \gamma_\eps(t,\mx)\leq C(|t-t_0|+|\mx-\mx_0|)\big(\delta_0(x_1)+r_\eps(t,\mx)\big)+ C.
\end{equation}
Here the uniform in $\eps$ constant $C$ depends on the Lipschitz constant of $F$, $G$, $\pa_u F$ and $\pa_u G$;
while  $(r_\eps)_\eps$ is a family of Radon measures
with total variation in $V(t_0,x_0)$ bounded by $M$, uniformly in $\eps$.
\end{lemma}
\begin{remark}\label{rem:Reps-1Dchoice}
  The analogous result can be stated for a discontinuous flux $\mff$ having a general Lipschitz jump manifold $\Sigma$.
  In this case, one has to assume that $R_\eps$ verify in $V(t_0,\mx_0)$ the viscous conservation law with flux coefficients $\mff^\pm$ from \eqref{eq:1sidedfluxes} frozen at the point $(t_0,\mx_0)$ on each side from the jump manifold $\Sigma$. Sometimes such profiles $R_\eps$ can be constructed, which may involve a
source term $S_\eps$ in the right-hand side that accounts for the curvature of $\Sigma$, see \cite{AKRnhm} and \S~\ref{ssec:CURVED-prelim} for more information. But in general, the construction of $R_\eps$ with some prescribed
piecewise constant behaviour, as $\eps\to 0$ (this is needed in Lemma~\ref{l2} below) is very delicate
when $\Sigma$ is not flat; while one-dimensional profiles $R_\eps$ with $S_\eps\equiv 0$
are readily constructed for a flat interface.
Therefore, in this paper we will stick to the basic choice $R_\eps=R(\frac{x_1}{\eps})$ for a flat boundary $\Sigma=\{(t,\mx)\,|\,x_1=0\}$, at a price of rectification arguments developed in Section~\ref{sec:CURVED} for reduction of curved manifolds $\Sigma$ to the flat case.
\end{remark}
\begin{proof}
First, we subtract \eqref{4} from \eqref{3}, in the resulting expression we add and subtract
the measure
$$
\Div_x \left(F(t,x,R_\eps)H(x_1)\! +\!G(t,x,R_\eps)H(-x_1)\right).
$$
We get
\begin{align*}
&\pa_t (u_\eps - R_\eps) +\Div_x \Big[
(F(t,\mx,u_\eps)-F(t,\mx,R_\eps))H(x_1)
\\&\qquad\qquad\qquad\qquad \ \ -
(G(t,\mx,u_\eps)-g(t,\mx,R_\eps))H(-x_1)\Big]\\
&+\Div_x \big( F(t,\mx,R_\eps)H(x_1)+ G(t,\mx,R_\eps)H(-x_1)\big)\\& -
\Div_x  \big( F(t_0,\mx_0,R_\eps) H(x_1)+
G(t_0,\mx_0,R_\eps)H(-x_1)\big)=\eps \Delta (u_\eps -R_\eps).
\end{align*}
This equality is understood in the distributional sense in
$V(t_0,\mx_0)$. Due to the $L^2(0,T;H^1_{loc}(\Omega))$ regularity
of $u_\eps,R_\eps$, proceeding by approximation we can multiply this
expression by  $\varphi\,{\rm sgn}_\alpha(u_\eps -R_\eps)$ where
${\rm sgn}_\alpha $ is a Lipschitz regularization of ${\rm sgn}$ and
$\varphi$ is a localizing test function. Classical chain-rule and
passage-to-the-limit in the regularization parameter $\alpha$
arguments apply (see, e.g., \cite{Otto-ellpb,Carrillo}). In this
way, we obtain in $\mathcal D'(V(t_0,\mx_0))$ the following
Kato-like inequality:
\begin{align*}
&\pa_t |u_\eps - R_\eps| +\Div_x \Big[ {\rm sgn}(u_\eps-R_\eps)\big(
(F(t,\mx,u_\eps)-F(t,\mx,R_\eps))H(x_1)\\&\qquad\qquad
+(G(t,\mx,u_\eps)-G(t,\mx,R_\eps))H(-x_1)\big)\Big]\;\leq\;|\omega_\eps|+\eps \Delta |u_\eps -R_\eps|,
\end{align*}

\vspace*{-16pt}
\begin{multline*}
\omega_\eps :=\Div_\mx \bigl[(F(t,\mx,R_\eps)- F(t_0,\mx_0,R_\eps)) H(x_1)\bigr]\\
 + \Div_\mx \bigl[(G(t,\mx,R_\eps)-G(t_0,\mx_0,R_\eps))  H(-x_1)\bigr].
\end{multline*}
 If we denote $\gamma_\eps = |\omega_\eps|$, we obtain \eqref{5}; it remains to estimate $\gamma_\eps$.

We get the bound \eqref{eq:bound-on-gamma} by computing $\omega_\eps$ explicitly. First, we estimate the jump term coming from the differentiation of $H(\pm x_1)$,
keeping in mind that $R_\eps$ is continuous across $\Sigma$ due to its $H^1_{loc}$ regularity in space.
Because $F(t,\mx,k)-F(t_0,\mx_0,k)$, $G(t,\mx,k)-G(t_0,\mx_0,k)$ are Lipschitz continuous functions
of $(t,\mx,k)$ taking value zero at $(t,\mx)=(t_0,\mx_0)$, the contribution of the jump term
is upper bounded by $C\,\text{dist}\big((t,\mx),(t_0,\mx_0)\big)\,\delta_0(x_1)$. Next, we focus on the terms coming from the differentiation of $F(t,\mx,R_\eps)-F(t_0,\mx_0,R_\eps)$ (the contributions of the analogous term with $G$ are estimated in the same way). We get the term $\Div_\mx F(t,\mx,r)|_{r=R_\eps(t,\mx)}$ bounded by the Lipschitz constant of $F$ and
the term
$$\nabla R_\eps(t,\mx)\,\cdot\,\big(\pa_u F(t,\mx,r)-\pa_u F(t_0,\mx_0,r)\big)|_{r=R_\eps(t,\mx)}$$
estimated by $r_\eps=|\nabla R_\eps|$ times the modulus of continuity of $\pa_u F$, $\pa_u G$ (that is how the distance between $(t,\mx)$ and $(t_0,\mx_0)$ enters the bound \eqref{eq:bound-on-gamma}). By assumption,  $\pa_u F$, $\pa_u G$
are Lipschitz continuous and  the integral of $r_\eps$ is bounded by the constant $M$. This leads to \eqref{eq:bound-on-gamma} and concludes the proof of the lemma.
\end{proof}

Next, we infer at the limit $\eps\to 0$ some information on traces of $u=\lim_{\eps\to 0} u_\eps$ that can be extracted from information available for the traces of $R= \lim_{\eps\to 0} R_\eps$:
\begin{lemma}
\label{l2} Let $V(t_0,\mx_0)$ be a neighbourhood of $(t_0,\mx_0)\in\Sigma$. In addition to the assumptions and notations of the
previous lemma, assume\\[3pt]
$\bullet$ $ u_\eps \to u$ while $R_\eps \to R$ a.e. in $V(t_0,\mx_0)$;\\[1pt]
$\bullet$ $u$ and $R$ admit strong left and right traces $u^-$ and
$R^-$, and $u^+$ and $R^+$, respectively, defined on $\Sigma\cap V(t_0,\mx_0)$;\\[1pt]
$\bullet$ and $(t_0,\mx_0)$ is the Lebesgue point of $u_{\pm}$ and of $R_{\pm}$ on $\Sigma$.
%

\smallskip
With the notation  $f$,$g$ introduced in~\eqref{eq:f-and-g},
 writing $u^\pm$ for $u^\pm(t_0,\mx_0)$ and $R^\pm$ for $R^\pm(t_0,\mx_0)$,
 there holds
\begin{equation}\label{7}
{\rm
sgn}(u^+-R^+)(f(u^+)-f(R^+))\,-\,{\rm
sgn}(u^--R^-)(g(u^-)-g(R^-))
\,\leq\, 0.
\end{equation}
\end{lemma}
\begin{proof} We first let $\eps\to 0$ in \eqref{5}, using the a.e. convergence and uniform boundedness of $u_\eps,R_\eps$ as well as \eqref{eq:bound-on-gamma}. We find that in $\mathcal D'(V(t_0,\mx_0))$, there holds
\begin{align}
\label{5-1} &\pa_t|u -R| +\Div_x \Big( {\rm sgn}(u \!-\!R)
((f(t,\mx,u)\!-\!f(t,\mx,R))H(x_1)\!
\\&\qquad\qquad\qquad\qquad\qquad\qquad \ \ +\!
(g(t,\mx,u)\!-\!g(t,\mx,R))H(-\!x_1)) \Big)\nonumber\\
& \leq C(|t-t_0|+|\mx-\mx_0|)\big(\delta_0(x_1)+r(t,\mx)\big)+ C. \nonumber
\end{align}
Here,  $r$ is a measure defined as a weak-* limit, along a subsequence, of the bounded sequence $(r_\eps)_\eps$ of Radon measures defined in Lemma~\ref{l1}. Now, proceeding by approximation with $C^\infty_c$ functions
%
 we test \eqref{5-1} with $\eta_h(t,\mx)=\varphi(t,\hat x_1) \mu_h(x_1)$,
where the function $\mu_h$ is defined by \eqref{mih} and
 $\varphi$ is regular, such that for $h$ small enough, $\eta_h$ is a Lipschitz function supported in $V(t_0,\mx_0)$.
 By letting $h\to 0$, using the existence of one-sided strong traces $u_\pm$, $R_\pm$ of $u,R$ on $\Sigma$ we get (with $\mx=(0,\hat x_1)$)
\begin{align*}
&\int_{\Sigma} \varphi(t,\mx) ({\rm
sgn}(u^+(t,\mx)-R^+(t,\mx))(F(t,\mx,u^+(t,\mx))-F(t,\mx,R^+(t,\mx)))dt
d\hat{x}_1\\&-\int_{\Sigma}\varphi(t,\mx) {\rm
sgn}(u^-(t,\mx)-R^-(t,\mx))(G(t,\mx,u^-(t,\mx))-G(t,\mx,R^-(t,\mx))))dt
d\hat{x}_1\\& \leq \int_{\Sigma}
C(|t-t_0|+|\mx-\mx_0|)(1+M)
\varphi(t,\mx)dt
d\hat{x}_1.
\end{align*}
Indeed, the contribution of the term $C$ to the right-hand side vanishes as $h\to 0$ and
the contribution of the measure $r(t,\mx)$ is estimated via the bound $M$ on its total variation.

  By taking for $\varphi$
an approximation of the Dirac measure concentrated at $(t_0,\mx_0)$
taking into account the fact that $(t_0,\mx_0)$ is the Lebesgue point of the functions $u_\pm$ and
$R_\pm$, we get the desired result.
\end{proof}

To continue, we need two more lemmas. The first
one provides a rough information on interface traces $u^\pm$ of $u$, in the spirit of \cite{kar3}.
The second one uses \eqref{7} of Lemma~\ref{l2} to describe more precisely the couples of possible interface traces, ensuring existence of a suitable value $p_u$ to be used in~\eqref{sep827}.
\begin{lemma}\label{lem:entropy-rough}
Assume that $u$ is an $L^1_{loc}$-limit (along a subsequence) of the family
$(u_\eps)_\eps$ of solutions to \eqref{3}, and it admits strong one-sided traces $u^\pm$ on $\Sigma$, as defined above.
With the notation  of Lemma~\ref{l2},
 the Rankine-Hugoniot condition $g(u^-)=f(u^+)$ holds, moreover, for all $k\in [a,b]$ there holds
\begin{equation}\label{rough-corollary}
 {\rm sgn}(u^+-k)(f(u^+)-f(k)) -{\rm sgn}(u^--k)(g(u^-)-g(k)) \leq |f(k)-g(k)|.
\end{equation}
\end{lemma}
\begin{proof}
The passage to the limit in the weak formulation \eqref{3} ensures that $u$ is a weak solution of \eqref{flat}, therefore the Rankine-Hugoniot condition~\eqref{R-H} of Lemma~\ref{lem:RH} holds.
Moreover, the following (rough) entropy inequality (which is precisely \eqref{eq:Panov-def} written for the case of equation~\eqref{flat}, see also~\eqref{eq:KRT-cond}) holds for every $k\in [a,b]$:
\begin{align}
\label{rough} &\pa_t|u-k|+\Div_\mx \Big({\rm sgn}(u-k)
((F(t,\mx,u)-F(t,\mx,k))H(x_1)\\\nonumber
&\qquad\qquad\qquad\qquad\qquad\qquad+(G(t,\mx,u)-G(t,\mx,k))H(-x_1))
\Big)\\&+{\rm sgn}(u-k) \Bigl(H(x_1)\Div_\mx
F(t,\mx,k)+H(-x_1)\Div_\mx G(t,\mx,k)\Bigr) \nonumber\\&
\leq |F_1(t,\mx,k)-G_1(t,\mx,k)|\delta_0(x_1), \nonumber
\end{align}
in the sense of distributions.
To prove~\eqref{rough}, following \cite{Kru} and \cite{pan_arma} it is enough to multiply \eqref{3} by a regularization of ${\rm sgn}(u_\eps-k)$ (cf. the proof of Lemma~\ref{l1}) with the rough estimation of the measure term charging $\Sigma$; then to let $\eps\to 0$ along the subsequence such that $u_\eps$ strongly converges
toward the function $u$.
Then, taking in~\eqref{rough} test functions $\mu_h(x)\xi(t,\mx)$ with $\xi\in \mathcal D(\Sigma)$ and $\mu_h$ given by \eqref{mih}, letting $h\to 0$ one finds that
\eqref{rough-corollary} holds in $\mathcal D'(\Sigma)$. The left-hand side of \eqref{rough-corollary}
being an $L^\infty(\Sigma)$ function, \eqref{rough-corollary} also holds pointwise, a.e. on $\Sigma$.
\end{proof}

\begin{lemma}\label{cruc}
Assume that $u$ is an $L^1_{loc}$-limit (along a subsequence) of the family
$(u_\eps)_\eps$ of solutions to \eqref{3}, and it admits strong one-sided traces $u^\pm$ on $\Sigma$, as defined above.  
Given $\mathcal H^d$-a.e point  $(t_0,\mx_0)\in \Sigma$, let $f,g$ and $u^\pm$ be defined as in Lemma~\ref{l2}.

\smallskip
Then there exists a measurable function $p_u: \Sigma \to [a,b]$ such that for $\mathcal H^d$-a.e.
 $(t_0,\mx_0)\in \Sigma$, the value $p_u=p_u(t_0,\mx_0)$ satisfies
\begin{multline}\label{sw-cond1}
{\rm sgn}(u^+-k)(f(u^+)-f(k)) -{\rm sgn}(u^--k)(g(u^-)-g(k))\\
 +{\rm sgn}(p_u-k)(f(k)-g(k))
\leq 0,
\end{multline}
 and $p_u$ lies in the (closed) interval with the
endpoints $u^+$ and $u^-$.
\end{lemma}
\begin{proof}
Fix $(t_0,\mx_0) \in \Sigma$, a Lebesgue point of traces $u^\pm$ defined on $\Sigma$. We will define pointwise $p_u(t_0,\mx_0)$; observe that the  construction given below is based on a case study of inequalities
 for functions $f=F(t_0,\mx_0,\cdot)$, $g=G(t_0,\mx_0,\cdot)$ that depend continuously on $(t_0,\mx_0)$; therefore it provides a measurable function $p_u$ defined on $\Sigma$.
 To simplify the analysis, assume
that $u^+\geq u^-$; the other case is symmetric. Remark that it is
enough to prove \eqref{sw-cond1} for $k\in [u^-,u^+]$ since for other values of $k$ the
result is evident.
 Indeed, if $k\notin (u^-,u^+)$ then with any choice of $p_u\in [u^-,u^+]$, inequality \eqref{sw-cond1} follows from the Rankine-Hugoniot relation $f(u^+)=g(u^-)$ of Lemma~\ref{lem:entropy-rough}.

We have the following cases:

\smallskip\noindent
$\bullet$ If the crossing condition \eqref{eq:crossing-cond} is satisfied for
the couple of functions $(f,g)$ on the interval $[u^-,u^+]$, then we define $p_u(t_0,\mx_0)$ (denoted $p_u$ in the sequel) by
\begin{equation*}
p_u:=
\begin{cases} u^+, &
f(k) \geq g(k) \ \ \forall k\in [u^-,u^+]\\
u^-, & f(k) \leq g(k) \ \ \forall k\in [u^-,u^+]\\
u^o, &  \text{if $u^o\in(u^-,u^+)$ is a crossing point of $f$ and $g$.}
\end{cases}
\end{equation*}
 Indeed, due to the crossing condition \eqref{eq:crossing-cond} the above choice yields
 $${\rm sgn}(p_u-k)(f(k)-g(k))=-|f(k)-g(k)|,$$
  which reduces \eqref{sw-cond1} to the already proved property \eqref{rough-corollary}.

\smallskip\noindent
$\bullet$  Assume next that
the couple of functions $(f,g)$ does not satisfy the crossing condition \eqref{eq:crossing-cond},
but there exists some intersection point (denoted again $u^o$) such that
\begin{align*}
&g(u^-) \leq g(k) \ \ k \in [u^-,u^o], \\
&f(u^+)\leq f(k), \ \
k \in [u^o,u^+].
\end{align*}
Then we can still proceed as above, taking $p_u=u^o$.

\smallskip\noindent
$\bullet$
Now, if none of the above possibilities holds, then there exists a crossing point
$u^o\in (u^-,u^+)$ of $f$ and $g$ such that there exist $k_-
\in (u^-,u^o)$, $k_+ \in (u^o,u^+)$
satisfying
\begin{align}
\label{e}
 g(u^-)> g(k_-) \ \
{\rm and} \ \
f(u^+)>f(k_+).
\end{align}
We will prove that this case is impossible, arguing by contradiction.
 First, notice that it holds:
\begin{equation}
\label{rel11}
g(u^-)=f(u^+)\leq
f(u^o)=g(u^o).
\end{equation} This relation is obtained by putting
$k=u^o$ in \eqref{rough-corollary} and using the Rankine-Hugoniot relation of Lemma~\ref{lem:entropy-rough}.
From \eqref{rel11} and the assumptions \eqref{e}, we see that we can
choose values $k^-\in [k_-,u^o)$ and $k^+\in (u^o,k_+]$ such that
\begin{equation}\label{nk1}
 g(k^-)=f(k^+)< f(u^+)=g(u^-) \leq g(u^o)=f(u^o);
\end{equation}
moreover, we will take
\begin{align*}
k^-=\max\Bigl\{\,k\in (u^-,u^o): \; g(k) \leq \max\{ g(k_-),f(k_+) \} \,\Bigr\}, \\
k^+=\min\Bigl\{\,k\in (u^o,u^+): \;
f(k) \leq \max\{
g(k_-),f(k_+) \} \,\Bigr\}.
\end{align*}

Now, using the procedure from the proof of \cite[Proposition
5.1]{AKR}, we argue that there exists a stationary solution to \eqref{4}
$R_\eps$ under the form of one-dimensional profile $R(\frac{x_1}{\eps})$ satisfying
\begin{equation}
\label{ic4} \lim_{z\to -\infty} R(z)=k^-, \ \ \lim_{z\to +\infty} R(z)=k^+, \
\ R(0)=u^o.
\end{equation} Indeed, $R$ is given by
$$
R(z)=
\begin{cases}
R^l(z), & z\leq 0\\
R^r(z), & z\geq 0
\end{cases},
$$
 where $R^l$ and $R^r$ are solutions to the ODEs
\begin{align*}
(R^l)'=g(k^-)-g(R^l), \
\ (R^r)'=f(R^r)-f(k^+)
\end{align*}
respectively,
with the initial data $R^l(0)=u^o=R^r(0)$. Remark that this does ensure
 $R^l(-\infty)=k^-$ and $R^r(+\infty)=k^+$. Indeed,
  according to the choice of $k^-$ and $k^+$, there holds
$g(R^l)-g(k^-)>0$, $R^l\in (k^-,u^o]$, and $f(k^+)-f(R^r)<0$,
$R^r\in [u^o, k^+)$. This actually means that $R^l$
will decrease from $0$ to $-\infty$ until it asymptotically reaches the stationary point $k^-$ of the corresponding ODE; while $R^r$ will decrease from $0$ to $\infty$ tending to $k^+$.

This construction provides $R \in  C(\R)\cap C^1(\R\setminus\{0\})$ and thus the corresponding rescaled profile $R_\eps$ belongs to $L^2_{loc}(\R^+;H^1_{loc}(\R^d))$ and it is readily checked that it represents a
solution to \eqref{4} in the sense of distributions.
Remark that such $(R_\eps)_\eps$ is a family of functions of uniformly
bounded variation; in view of the explicit formula \eqref{ic4},
$R_\eps$ converges to $R(x_1)=k^-$ for $x_1<0$, and to $R(x_1)=k^+$ for $x_1>0$. Therefore  the family of functions $(R_\eps)_\eps$ satisfies conditions of Lemmas~\ref{l1},\ref{l2} with $S_\eps=0$, the corresponding traces $R^\pm$ of $R=\lim_{\eps\to 0} R_\eps$ being equal to $k^\pm$. The family $(u_\eps)_\eps$ satisfies as well the assumptions of these lemmas, therefore we conclude from
\eqref{7} (recalling the meaning of notation for $f,g$, $u^\pm$, $R^\pm$)
that the following relation holds
\begin{equation*}
{\rm
sgn}(u^+ -k^+)(f(u^+)-f(k^+)) -{\rm sgn}(u^--k^-)(g(u^-)-g(k^-)) \leq 0.
\end{equation*}
Since $u^-\leq k^- \leq k^+ \leq u^+ $, this reduces to
$$
f(u^+)+g(u^-)\leq f(k^+)+g(k^-),
$$
which contradicts \eqref{nk1} implying that case \eqref{e} is not possible.
This concludes the proof.
\end{proof}

Now, it is easy to prove existence of the entropy admissible
solution to \eqref{flat}, \eqref{flati}.

\begin{theorem}
\label{thm5} Let \eqref{gnl},\eqref{eq:confinement},\eqref{eq:Lip-flux} and \eqref{hyp:pa_u F,G continuous in t,x} be fulfilled.
 Then, there exists a solution $u$ to \eqref{flat}, \eqref{flati} in the sense of Definition \ref{def-adm}.
\end{theorem}

\begin{proof}
 Existence of weak solutions to~\eqref{3}  in $L^2_{loc}(\R^+;H^1_{loc}(\R^d))$ (as required in Lemma~\ref{l1}) with a given initial datum $u_0$  can justified with rather classical arguments. Let us give a very brief sketch of stages of one among many possible construction arguments.

 \smallskip\noindent
 $\bullet$ One can start with $L^2(\R^d)$ data $u_0$ and construct solutions $u_\eps$ in the energy space $L^2(0,T;H^1(\R^d))$ with Galerkin approximations, see, e.g., \cite{Otto-ellpb} and references therein,
 see also \cite[Sect.\,6.2]{AKR}.

%

  \smallskip\noindent
 $\bullet$ Further, assumption \eqref{eq:confinement}  along with the comparison principle for $u_\eps$ (that can be proved as in \cite{Otto-ellpb,Carrillo,Vallet,BendahmaneKarlsen}) ensures the confinement
 property $a\leq u_\eps\leq b$, which yields the uniform $L^\infty$ bound on $(u_\eps)_\eps$.

  \smallskip\noindent
 $\bullet$ Using, by approximation, the test function $\exp(-|x|)u_\eps$,
 one finds  estimates on $u_\eps$ in the space $L^2(0,T;H^1(\R^d,w))$ for all $T>0$, where $H^1(\R^d,w)$ is the weighted $H^1$ space with the norm defined by $\|u\|^2_{1,w}:=\int_{\R^d} \exp(-|x|)\bigl(|u|^2+|\nabla u|^2\bigr)$;
  as soon as $u_0$ remains $[a,b]$-valued, these estimates do not depend on $\|u_0\|_{L^2(\R^d)}$.

  \smallskip\noindent$\bullet$ With a.e. approximation by $[a,b]$-valued $L^2(\R^d)$ functions, using space and time compactness that follows from the above uniform estimates and the variational formulation of \eqref{3} in $\bigl(L^2(0,T;H^1(\R^d,w))\bigr)^*$ (cf., e.g., \cite{Otto-ellpb}),
  one can extend the solution construction to any $[a,b]$-valued $u_0$ and get solutions $u_\eps\in L^2(0,T;H^1(\R^d,w))$. Clearly, these solutions lie in $L^2_{loc}(\R^+;H^1_{loc}(\R^d))$.

 \smallskip
  Then, \eqref{gnl} guarantees existence of an accumulation point $u$ of $(u_\eps)_\eps$, as $\eps\to 0$, thanks to strong precompactness theorems (\cite{pan_arma}, see also~\cite{MA}).

 By construction, the result of the previous lemma can be applied to $u$.
It is not difficult to check that the function $p_u$ defined in the
previous lemma and the strong limit $u$ of $(u_\eps)_\eps$ the vanishing viscosity
approximations of \eqref{3} with fixed initial datum $u_0$ satisfy conditions of
Definition \ref{def-adm}.

Indeed,
given $k\in [a,b]$, denote by $L_k$ the element of the dual of $C^1_c(\R^+\times\R^d)$ defined by the left-hand side of \eqref{sep827}.
For any nonnegative compactly supported function
 $\varphi \in \mathcal C^1(\R^+\times \R^d)$, using $\mu_h$ defined by \eqref{mih} we have:
\begin{equation}\label{eq:splitting}
 L_k(t,x) \bigl(\varphi\bigr)= L_k\bigl(
(1-\mu_h)\varphi\bigr) + L_k\bigl(\mu_h\varphi\bigr).
\end{equation}
The first summand in the right-hand side of \eqref{eq:splitting} is less than or equal to zero
due to \eqref{rough}, since the truncation function $(1-\mu_h)$ is supported on a subset of $\R^+\times\R^d$ where the left-hand sides of \eqref{rough} and \eqref{sep827} coincide. 
 Further, due to existence of traces $u^\pm$
  and to the continuity of the maps \eqref{eq:continuity-of-KrFluxes}
 the limit as $h\to 0$ of the second summand in the right-hand side of \eqref{eq:splitting} can be written as
 $$
 \int_{\Sigma} \Delta_u(t_0,\mx_0)\varphi(t_0,\mx_0)\,dt_0 d\hat{x}_{10}
 $$
 which is non-positive, since the quantity
 \begin{multline*}
   \Delta_u(t_0,\mx_0):={\rm sgn}(u^+(t_0,\mx_0)-k)(f(u^+(t_0,\mx_0))-f(k))\\ -{\rm sgn}(u^-(t_0,\mx_0)-k)(g(u^-(t_0,\mx_0))-g(k))
 +{\rm sgn}(p_u(t_0,\mx_0)-k)(f(k)-g(k))
 \end{multline*}
is non-positive due to inequality~\eqref{7} proved in Lemma~\ref{cruc}.
Thus, as $h\to 0$ we find that $L_k$ is a non-positive distribution, which proves~\eqref{sep827}.
\end{proof}

\subsection{Remark to the definition of admissibility}\label{ssec:FLAT-pu-specified}

Let us underline that a shortcoming of the admissibility concept of Definition~\ref{def-adm}
is that the singular values given by a function $p_u:\R^+\to [a,b]$ are not explicitly determined on the set $\Sigma$, see Remark~\ref{rem:micromodel}(iv).
To weaken this indetermination, let us show that the function $p_u$ on $\Sigma$ can be
chosen so that it takes values only from the set of the crossing
points of the fluxes $f,g$ defined in~\eqref{eq:f-and-g} (observe that $a,b$ always belong to the set of crossing points, due to assumption~\eqref{eq:confinement}).
 Indeed, we have
\begin{proposition}
\label{ekv} Assume that the function $u:\R^+\times \R^d \to [a,b]$
represents an entropy admissible solution to \eqref{flat}, i.e., \eqref{sep827} holds with some $\mathcal H^d$-measurable function $p_u$ on $\Sigma$.
 Then
there exists another $\mathcal H^d$-measurable function $p_u:\Sigma\to [a,b]$ such that for $\mathcal H^d$-a.e. point
$(t_0,\mx_0)\in \Sigma$, one has
 $$p_u(t_0,\mx_0)\;\in\;  {\cal C}(t_0,\mx_0):=\Bigl\{k\in [a,b]: \; F_1(t_0,\mx_0,k)=G_1(t_0,\mx_0,k)\Bigr\},$$
such that \eqref{sep827} is still satisfied with this new function $p_u$.
\end{proposition}

\begin{proof} Notice that it is enough to prove that there exists required function $p_u$ such that \eqref{sw-cond1} holds for the
right and left traces $u^+$ and $u^-$, respectively, of the function
$u$. Since $u^-$ and $u^+$ must satisfy \eqref{R-H}, we have the
following possible situations ($\mathcal H^d$-a.e. point of $(t_0,\mx_0)\in \Sigma$ being fixed in the arguments that follow), classified as in the case study \eqref{i1}--\eqref{vi2}:

\medskip\noindent
{\bf Case 1:} $ \ \ u^+\leq u^-$.

\smallskip\noindent
$\bullet$ $u^+\leq u^-\leq p_u$.\\
From \eqref{i1} and \eqref{ii1} we see that instead of $p_u$ we can take
$u^o \in {\mathcal C}$ such that it is the minimal
element from ${\mathcal C}$ satisfying $u^o \geq p_u$.

\smallskip\noindent
$\bullet$ $u^+\leq p_u \leq u^-$.\\
From \eqref{iii1} and \eqref{iv1} and the Rankine-Hugoniot
conditions, we see that there exists $u^o\in {\mathcal
C}$ such that $u^+\leq u^o \leq u^-$. Such $u^o$ can
be taken instead of the given $p_u$.

\smallskip\noindent
$\bullet$ $p_u \leq u^+\leq u^-$.\\
From \eqref{v1}  and \eqref{vi1} we see that instead of $p_u$ we can take
$u^o\in {\mathcal C}$ such that it is the maximal
element from ${\mathcal C}$ satisfying $u^o \leq p_u$.

\medskip\noindent
{\bf Case 2:} $ \ \ u^-\leq u^+$.

 \smallskip\noindent
$\bullet$ $u^-\leq u^+\leq p_u$.\\
From \eqref{i2} and \eqref{ii2} we see that instead of $p_u$ we can take
$u^o\in {\mathcal C}$ such that it is the minimal
element from ${\mathcal C}$ satisfying $u^o \geq p_u$.

\smallskip\noindent
$\bullet$ $u^-\leq p_u \leq u^+$.\\
From  \eqref{iii2} and \eqref{iv2}  and the Rankine-Hugoniot
conditions, we see that there exists $u^o\in {\mathcal
C}$ such that $u^+\leq u^o \leq u^-$. Such $u^o$ can
be taken instead of the given $p_u$.

\smallskip\noindent
$\bullet$ $p_u \leq u^-\leq u^+$.\\
From  \eqref{v2} and \eqref{vi2} we see that instead of $p_u$ we can take
$u^o\in {\mathcal C}$ such that it is the maximal
element from ${\mathcal C}$ satisfying $u^o \leq p_u$.

\medskip
This concludes the proof, the value $p_u$ being re-defined in each case by $p_u:=u^o \in \mathcal C$.
Observe that the measurability of the so re-defined $p_u$ follows from the continuity of $F_1$,$G_1$ as functions
of $(t_0,\mx_0)\in\Sigma$ and from measurability of the traces $u^\pm$ and of the original function $p_u$ defined on $\Sigma$.
\end{proof}

\section{Piecewise regular discontinuous-flux Cauchy problems}\label{sec:CURVED}

Clearly, the flat case of the previous section is a model case. More realistic problems of the kind \eqref{eq:cons-law} can present flux discontinuities along multiple curved, possibly intersecting hypersurfaces: one example is the case of flows in a homogeneous porous medium with inclusions of another homogeneous porous medium.
In this section, we justify well-posedness of
\eqref{eq:cons-law}, \eqref{flati} in the case where the
discontinuity interfaces $\Sigma$ are not necessarily flat, with two different approaches to construction of solutions under \emph{ad hoc} technical restrictions on the singularities of the flux $\mff$.

\subsection{Guidelines for extension of well-posedness results}\label{ssec:CURVED-prelim}

First, let us explain why the uniqueness proof is generalized in a straightforward way, and what are the difficulties of the existence proof and the ideas to overcome these difficulties.

We start with state-of-the-art observations.
Approaches to uniqueness and existence in the case of multiple, non-flat and possibly crossing
jump manifolds $\Sigma$ of the flux $\mff$ have already been developed (in the setting of the Karlsen-Risebro-Towers definition, \cite{kar3}; and in the setting of $\mathcal G_{VV}$-entropy solutions, \cite{AKRnhm}). In a part,
they can be exploited with our new definition \eqref{eq:Panov-def-extended},\eqref{eq:p-function}.

 \medskip
 \noindent$\bullet$ Both the definition of admissible solution and the uniqueness proof of the flat case are easily extended to the setting of piecewise Lipschitz regular, genuinely nonlinear fluxes
with a locally finite number of jump discontinuities along Lipschitz
hypersurfaces (see \cite{AKRnhm}, see also \cite{kar3}).

For instance, following the approach of \S~\ref{ssec:FLAT-GVV-uniq},
one readily extends the definitions to the case of a curved interface $\Sigma$: it is enough to replace
the couple $(f,g)=(F_1,G_1)$ in Definition~\ref{gvv-adm} by the couple of normal flux functions
with coefficients frozen at $(t_0,\mx_0)$:
 \begin{equation*}
\begin{array}{l}
 \text{$\mff:u\mapsto F(t_0,\mx_0,u)\cdot\nu(t_0,\mx_0)$
(flux in the direction of $\nu(t_0,\mx_0)$),}\\
 \text{$g:u\mapsto G(t_0,\mx_0,u)\cdot\nu(t_0,\mx_0)$ (flux in the opposite direction)},
\end{array}
 \end{equation*} where $\nu(t_0,\mx_0)$ is some fixed unit normal vector to $\Sigma$ at the point $\mx_0$, and traces $u^\pm(t_0,\mx_0)$ are taken according to this fixed orientation of $\nu(t_0,\mx_0)$. As in the flat case, being solution in the sense of Definition~\ref{gvv-adm} is equivalent to being admissible solution in the sense of \eqref{eq:Panov-def-extended},\eqref{eq:p-function}. Indeed, the arguments that ensure constraints on the couple $(u^-,u^+)(t_0,\mx_0)$ starting from \eqref{eq:Panov-def-extended},\eqref{eq:p-function} reduce to a pointwise discussion.

 For the uniqueness proof, first observe that the analogue of the Kato inequality \eqref{kato-per}
is easily obtained for a locally finite union $\Sigma$ of Lipschitz interfaces. Then the constraints
on $(u^-,u^+)(t_0,\mx_0)$ derived from \eqref{eq:Panov-def-extended} (constraints that amount to the fact that $(u^-,u^+)(t_0,\mx_0)$
belongs to the vanishing viscosity germ associated with the point $(t_0,\mx_0)$ of $\Sigma$)
give the desired sign to the quantity $\Delta(u^\pm,v^\pm)(t_0,\mx_0)$ in the analogue of inequality \eqref{kato-per}. This results in the $L^1$-contraction principle (with an exponential growth term,
if a non-zero, Lipschitz continuous in $u$ source term $S$ is present in equation~\eqref{eq:cons-law}) and to uniqueness of an admissible solution in the sense \eqref{eq:Panov-def-extended},\eqref{eq:p-function}.

\smallskip
\noindent$\bullet$ At the same time, in order to justify the existence of solutions in the sense \eqref{eq:Panov-def-extended},\eqref{eq:p-function},
 unlike for the Karlsen-Risebro-Towers definition (\cite{kar3}) we cannot hope for a simple existence proof based on straightforward passage to the limit from approximate solutions: nothing seems to ensure strong compactness of singular values $(p_{u_\eps}(\cdot))_\eps$ on $\Sigma$.
  An indirect existence proof analogous to that of the flat case in~\S~\ref{ssec:FLAT-exist} is known (\cite{AKRnhm}, see also Remark~\ref{rem:Reps-1Dchoice})
  under the additional local regularity of jump manifolds: for instance, piecewise convexity or concavity of these manifolds is enough. The construction was proposed in \cite{AKRnhm} where it was presented under the simplifying assumption that the flux is piecewise constant in $(t,\mx)$ and the jump manifold is given by the graph of a $t$-independent $C^2$ function: $x_1=\phi(\hat x_1)$.
It is straightforward to extend this method to $(t,\mx)$-dependent genuinely nonlinear fluxes with jump manifolds given as portions $\Sigma_i$ of graphs of the form $x_1=\phi_i(t,\hat x_1)$ (up to relabeling the axes), where functions $\phi_i$ are Lipschitz continuous in $(t,\mx)$ and such that, additionally, $\Delta_{\hat x_1} \phi_i$ is a Radon measure on $\R^+\times \R^{d-1}$. We refer in particular to \cite[estimate (28)]{AKRnhm} for the details of the computations.

\smallskip
\noindent$\bullet$
To sum up, uniqueness of solutions in the sense \eqref{eq:Panov-def-extended},\eqref{eq:p-function} for \eqref{eq:cons-law},\eqref{eq:confinement} can be proved for the case of piecewise Lipschitz $\mff$ with locally finite number of Lipschitz components in the jump manifolds $\Sigma$. However, existence in this sense requires unnatural assumptions, if we follow the method of~\cite{AKRnhm}; while direct existence arguments used, e.g., in \cite{kar3}, are not applicable.

\medskip
In this paper, we want to drop the unnatural restriction on $\Delta_{\hat x_1} \phi_i$ in the existence proof
based on the approach of \cite{AKRnhm} and \S~\ref{ssec:FLAT-exist}.
 The method mainly consists in comparison of the solution with limits of appropriate viscosity profiles.
  Yet we do not attempt to construct approximate viscosity profiles (cf. Remark~\ref{rem:Reps-1Dchoice}) satisfying the original equation up to a measure source term not charging the interface (such construction is successful for the flat case or under the above mentioned restrictions on  $\Delta_{\hat x_1} \phi_i$,  but it is very delicate in the general case).
Instead,  we rely upon the invariance of the notion of solution - in the sense \eqref{eq:Panov-def-extended},\eqref{eq:p-function} -  under Lipschitz changes of coordinates
 in equation \eqref{eq:cons-law}. This invariance permits us to exploit the possibility to rectify $\Sigma$, at least locally, with a change of variables.
 Actually we propose two methods for constructing solutions: the one in \S~\ref{ssec:CURVED-AlmostRect} (given with the essential details and examples)  and the one of \S~\ref{ssec:CURVED-LocalDiffusion} (this method is more general but the details of construction are quite tedious, therefore they are only sketched).

%
%
As a matter of fact, in both cases the construction of solutions via a change of variables
amounts to a specific ``adapted'' viscosity approximation of the original equation \eqref{eq:cons-law}.
Such approximation is considerably more involved than the isotropic homogeneous viscosity $\eps \Delta u_\eps$ used
in most of the previous works on the subject (see, e.g., in \cite{kar3,pan_arma,Die3,AKRnhm}).
The use of an involved viscosity regularization operator need not be seen as a drawback of the method: indeed, we put forward the following remark.
\begin{remark}\label{rem:viscosity-bis}  As observed in Remark~\ref{rem:micromodel}(ii),
an important feature of the vanishing viscosity approximation is to ensure the global
continuity of the approximate solutions (this contrasts with the properties of the vanishing capillarity regularizations, see, e.g., \cite{Kaas,AndrCances}). Regularizing the strongly heterogeneous
first-order model~\eqref{eq:cons-law} with introduction of the
isotropic homogeneous diffusion $\eps \Delta u$ may look un-realistic.
 Heuristically, along with
this basic vanishing heterogeneous diffusion operator one should
consider, for the existence proof, the possibility to approximate
solutions via more involved heterogeneous diffusion operators that
still ensure global continuity of $u_\eps$.
\end{remark}

\subsection{Existence for almost rectifiable sets of singularities}\label{ssec:CURVED-AlmostRect}

The proof developed here relies on three ingredients.
 First, we have already known how to construct solutions (where the main difficulty is to construct $p_u$)
 for the case of a flat portion of the interface; moreover, the arguments put forward in
 Lemmas~\ref{l1}--\ref{cruc} are local arguments that work in a neighbourhood of $\mathcal H^d$-a.e. interface point,
 provided the corresponding portion of the interface is flat.
 Second, we point out the invariance of the notion of admissible solution
 with respect to a family of (possibly singular) changes of variables.
 Third, it is possible to deal with interface-rectifying changes of variables involving lower-dimensional singularities.

 The combination of the two latter ideas leads to the definition of \emph{almost rectifiability} for jump manifolds,
illustrated by several examples that show that this notion is in fact rather general.
\begin{definition}\label{def:almost-rec}
A union of hypersurfaces $\Sigma$ in $\R^+\times\R^d$ is called almost rectifiable
if there exists a direction  $\vec e_1\in \R^d$ (that, up to a
rotation of coordinate axes, can be assumed to be the direction of
the first canonical basis vector, i.e., $\pa_{\vec e_1}=\pa_{x_1}$)
and a closed set $\gamma\subset\Sigma$ consisting of a locally
finite number of disjoint Lipschitz manifolds
of codimension in $\R^{d+1}$ greater than one, satisfying the following $(\Sigma\setminus\gamma)$-rectification property:\\[3pt]
There exists a continuous surjective map
$$
\Phi:\R^+\times\R^{d}\longrightarrow \R^+\times\R^d,\;\;
\Phi:(\tilde t, \tilde x_1, \widetilde{\hat x}_1) \mapsto
(t,x_1,\hat x_1)
$$
that only acts in the direction $\vec e_1$, i.e., it has the form
\begin{equation}\label{eq:changeofvariables}
   t=\tilde t,\;\; \hat x_1=\widetilde{\hat x}_1,\;\;  x_1=\phi(\tilde{t},\tilde{\mx})
\end{equation}
with some scalar function $\phi$ on $R^+\times \R^d$,
and such that $\Phi$ fulfills the following properties:

\smallskip\noindent
(a) There exists a closed subset $\Gamma$ of $\R^+\times\R^{d}$ such that
   \begin{enumerate}
     \item One has $\Phi(\Gamma)=\gamma$, moreover, the restriction of $\Phi$
     to $\tilde\Omega:=(\R^+\times\R^{d})\setminus\Gamma$
     is injective;
     \item the map $\Phi$ is Locally Lipschitz continuous on $\tilde \Omega$ and
   the inverse of $\Phi$, defined on  $\Omega:=(\R^+\times\R^{d})\setminus\gamma$,
  is locally Lipschitz continuous as well;
  \item  in addition, the second-order derivatives $\pa^2_{\tilde{t}\tilde{x}_1} \phi$ and $\pa^2_{\tilde{x}_j\tilde{x}_1} \phi$, $2\leq j\leq d$,  are locally bounded on $\tilde{\Omega}$;
   \end{enumerate}

\noindent
  (b) There exists a set $\tilde{\Sigma}\subset\tilde\Omega$ consisting of flat
  manifolds orthogonal to the direction $\tilde{x}_1$ such that $\Phi(\tilde{\Sigma})=\Sigma\setminus\gamma$.
\end{definition}

%
%

\begin{figure}[htp]
\begin{center}
  \includegraphics[width=4in]{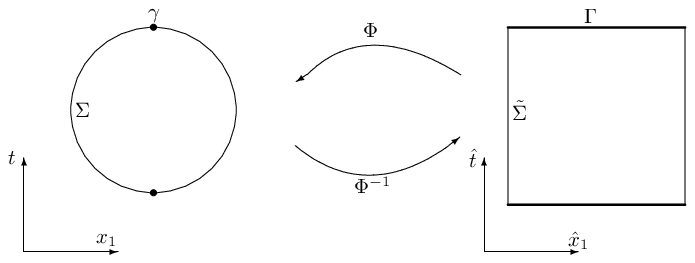}
  \caption{Change of variables transforming square into circle}
  \end{center}
\label{fig:CircleToSquare}
\end{figure}

Heuristically, the almost-rectifiability property of $\Sigma$ means that, up to a union of lower-dimensional
sub-manifolds $\gamma$ of the union $\Sigma$ of hypersurfaces of $\R^{d+1}$, a change of variables
rectifies $\Sigma$. To do this, one may need  $\Phi^{-1}$ to be multivalued.
 For instance,
a circle can be transformed into a square by rectifying the two open half-circles,
but then the two poles of the circle have to be stretched into segments (see Figure~\ref{fig:CircleToSquare}).
\begin{remark}\label{rem:changeofvar}
To construct such transformation in practice, having chosen coordinates in such a way that $\vec e_1$ is transversal to $\Sigma$ in  a.e. point with respect to the $d$-dimensional Hausdorff measure on $\Sigma$,
the map $\phi_{t,\tilde{\hat
{x}}_1}(\cdot):=\phi(t,\cdot,\tilde{\hat
{x}}_1)$ from $\R$ to $\R$ can taken to be a non-strictly increasing, piecewise polynomial function;
 moreover, whenever $\Sigma$ is a union of Lipschitz hypersurfaces, one can ask for a Lipschitz dependence of the function $\phi_{t,\tilde{\hat {x}}_1}(\cdot)$ on the parameters $t,\tilde{\hat
{x}}_1$ (such construction leads to the desired bounds on the mixed second derivatives $\pa^2_{\tilde{t}\tilde{x}_1} \phi$ and $\pa^2_{\tilde{x}_j\tilde{x}_1} \phi$).

The simplest explicit example is the case where $\Sigma:=\{(t,x_1,\hat x_1)\,|\, x_1=\phi_0(t,\hat x_1)\}$
in which case it is enough to set $x_1=\tilde x_1 + \phi_0(t,\widetilde{\hat x}_1)$. The properties (a)(1)-(3) and (b)
are readily checked in this case. Further examples are presented in \S~\ref{ssec:CURVED-rectifiable-examples}.
\end{remark}
In the sequel, whenever convenient we will denote $(\tilde x_1, \tilde{\hat x}_1)$ by $\my$ and we will not distinguish $\tilde t$ and $t$ since they are equal; to summarize the definition, we can write
 $$
 \begin{array}{l}
   (t,\mx)=\Phi(\tilde t,\tilde x_1, \tilde{\hat x}_1)=\Phi(t,\my),\\
    \Phi \; \text{is a locally bi-Lipschitz bijection from $\tilde\Omega=(\R^+\times\R^{d})\setminus\Gamma$}
    \\
    \text{to $\Omega=(\R^+\times\R^{d})\setminus\gamma$ such that $\Phi^{-1}$ rectifies $\Sigma\setminus\gamma$},
 \end{array}$$
 with specific mixed second derivatives of $\Phi$ that are locally bounded on $\tilde\Omega$.

\smallskip
 In the sequel, we assume that the flux $\mff$ in
\eqref{eq:cons-law} has its jump manifolds contained in an almost
rectifiable set $\Sigma$ and it is Lipschitz continuous in all variables on each connected component
of the complementary of $\Sigma$.
Notice that one can drop the assumption $\gamma\subset\Sigma$
whenever this is convenient for construction of the change of
variables $\Phi$: indeed, one can always extend $\Sigma$ by
including ``ghost'' interfaces across which the flux $(t,\mx)\mapsto
\mff(t,\mx,k)$ has jump zero, for all $k$.

\smallskip
Using the change of variables $\Phi$ described in Definition~\ref{def:almost-rec},
 we transform equation \eqref{eq:cons-law} set up in $\R^+\times\R^d$ into an analogous equation set up in the domain
$\tilde\Omega=\Phi^{-1}(\Omega)$, where $\Omega$ is the set $(\R^+\times\R^d)\setminus \gamma$
which complementary $\gamma$ is of codimension $2$ or more. Indeed, observe that considering the first-order
equation \eqref{eq:cons-law} (in a weak or in an entropy sense) in the whole space $\R^+\times\R^d$  is equivalent to considering it
 (in the same sense) restricted to $\Omega$: cf. the proof of Proposition~\ref{prop:equiv-changeofvar} below.
Then, by an explicit calculation we see that under the change of variables $(t,\mx)=\Phi(t,\my)$
equation \eqref{eq:cons-law}
 becomes the  equation on $\tilde \Omega$ under the analogous conservative form
\begin{equation}
\label{general_app4}
 \pa_t \tilde{u} +\text{div}_{\my}
\tilde \mff (t,\my,\tilde{u})=\tilde S(t,\my,\tilde{u}).
\end{equation}
 Here
 $$
 \tilde \mff(t,\my,k)=\mathfrak D(t,\my)\,f(\Phi(t,\my),k)
 + k \pa_t \phi(t,\my)\vec e_1
 $$
 ($\mathfrak D(t,\cdot)$ being the $d\times d$ Jacobian matrix of the change of variables $\mx$ into $\my$ given by $\mx=\Phi(t,\my)$)
 and $\tilde S$ is computed using $\mff,S$ and second-order derivatives of $\phi$ that appear in assumption (a)(3) of Definition~\ref{def:almost-rec}.
For the Ansatz \eqref{eq:changeofvariables} and under the assumption (a) of Definition~\ref{def:almost-rec},
the scalar function $\tilde S(t,\my,k)$
defined on $\tilde \Omega\times [a,b]$
is measurable with respect to $(t,\my)$ and locally Lipschitz continuous in $k$.
Further,
  $\tilde \mff=\tilde \mff(t,\my,k)$ is a piecewise continuous vector-function
defined on $\tilde\Omega\times[a,b]$ which is discontinuous with respect
to $(t,\my)$ along the union of flat manifolds $\tilde \Sigma$
contained in $\tilde \Omega$ (see (b) of
Definition~\ref{def:almost-rec}),
and which is Lipschitz with respect to $(t,\my,k)$, locally in $(t,\my)\in\tilde \Omega$ and globally in $k\in[a,b]$.
Similarly, if one imposes the $(t,\mx)$-Lipschitz property \eqref{hyp:pa_u F,G continuous in t,x} of the derivatives $\pa_u \mff$ in each component of  $(\R^+\times\R^d)\setminus\Sigma$, this property persists for the flux $\tilde \mff$ in variables $(t,\my)$, in each component of  $\tilde\Omega \setminus \tilde \Sigma$, locally in $\Omega$ (i.e., the Lipschitz constant may blow up as $(t,\mx)$ approaches $\Gamma$).

\smallskip
Furthermore, the genuine nonlinearity property \eqref{gnl} is inherited by the flux $\tilde \mff$ on $\tilde\Omega\times [a,b]$ provided it is satisfied by $\mff$ on $(\R^+\times\R^{d})\times [a,b]$.
Finally, property \eqref{eq:confinement} for $\mff,S$ ensures that for equation \eqref{general_app4},
\begin{equation}\label{eq:a,b=sub,supersolution}
\text{the constant $a$ (resp, $b$) is a weak sub- (resp, super-) solution}.
\end{equation}

\smallskip
 Because equation \eqref{general_app4} has to be considered only in $\tilde \Omega=(\R^+\times\R^d)\setminus \Gamma$ (recall that $\Omega\cap \{(t,\mx)\,|\,t>0\}$ is an open subset of $\R^+\times \R^d$ which can have holes or cracks, see Figure~\ref{fig:CircleToSquare} and examples below), i.e., with test functions that belong to $\mathcal D(\tilde \Omega)$, we propose the following definition
 of which Definition~\ref{def-adm} is a particular case corresponding to $\Gamma=\emptyset$ and $\Sigma=\{(t,\mx)\,|\,x_1=0\}$.
\begin{definition}\label{def-mman}
 We say that a function $\tilde u$ taking values in
$[a,b]$ is an admissible solution to \eqref{general_app4} with
initial datum $\tilde u_0$ if $\tilde u$ it satisfies the analogue
of \eqref{eq:Panov-def-extended} within $\tilde \Omega$:
\begin{multline}\label{eq:Panov-def-extended-bis}
    \forall k\in\R\;\; |\tilde u\!-\!k|_t \; +\; \text{div}_\my\Bigl(\text{sgn}(\tilde u\!-\!k)(\tilde \mff(t,\my,\tilde u)\!-\!\tilde \mff(t,\my,k)) \Bigr)\\
    \leq\; - \text{sgn}(p_{\tilde u}\!-\!k)\text{div}_\my \mff(t,\my,k)\;\; \text{in $\mathcal D'(\tilde \Omega)$}
  \end{multline}
  with $\tilde u(0,\cdot)=\tilde u_0$,
where
\begin{equation}\label{eq:p-function-bis}
\text{$p_u$ is some Borel function that coincides with $\tilde{u}$
a.e. on $\tilde\Omega$.}
\end{equation}
\end{definition}
As in Definition~\ref{def-adm}, the initial datum entering the
entropy inequalities can be included in $C(\R^+;L^1_{loc}(\R^d))$
sense. The following easy observation follows:
\begin{proposition}\label{prop:equiv-changeofvar}
Assume that the flux $\mff$ in \eqref{eq:cons-law} is piecewise Lipschitz, discontinuous
in $(t,x)$ across a set $\Sigma$  almost rectifiable in the sense of
Definition~\ref{def:almost-rec}.\\[3pt]
Given a measurable function $u_0$ on $\R^d$ taking values in
$[a,b]$, consider initial condition $\tilde u_0$ defined on $\tilde
\Omega_0=\partial\tilde\Omega\cap\{(t,\my)\,|\,t=0\}$ by $\tilde u_0(\my):=u_0(\mx)$ where $\my$ is such that
$\mx=\Phi(0,\my)$. Let $\tilde \mff,\tilde S$ be defined so that
\eqref{eq:cons-law} transforms into \eqref{general_app4}
under the change of variables $\Phi$.\\[3pt]
A function $\tilde u$ on $\tilde \Omega$ is an admissible solution
of \eqref{general_app4} with initial datum $\tilde u_0$ in the sense
of Definition~\ref{def-mman} if and only if the function
$u$ on $\R^+\times\R^d$ such that $u(t,\mx)=\tilde u(t,\my)$ with $\mx=\Phi(t,\my)$ is an admissible solution
of \eqref{eq:cons-law}  in the sense
 \eqref{eq:Panov-def-extended},\eqref{eq:p-function} with the initial datum $u_0$.
\end{proposition}
\begin{proof}
 Since $u$, $\mff(\cdot,u(\cdot))$ and $S(\cdot,u(\cdot))$ belong to $L^\infty(\R^+\times\R^d)$ and $\tilde u$, $\tilde \mff(\cdot,\tilde u(\cdot))$ and $\tilde S(\cdot,\tilde u(\cdot))$ belong at least to $L^\infty_{loc}(\tilde\Omega)$,
the spaces $\mathcal D(\Omega)$ and $\mathcal D(\tilde\Omega)$ for test functions can be replaced by the spaces $W^{1,1}_c$
of locally Lipschitz, compactly supported in $\Omega$ and in $\tilde\Omega$, respectively, test functions.
The latter spaces are transported the one onto the other under the change of variables $\Phi$ that maps $\tilde\Omega$ onto $\Omega$.

Therefore, first, the claim of the proposition holds true when
the inequalities \eqref{eq:Panov-def-extended},\eqref{eq:p-function} for $u$ are restricted to $\mathcal
D'(\Omega)$ test functions. Indeed, the change of variables transforms one inequality into the other,
with the transformed test function. In particular, under the change $\Phi$
the values $p_u$ on $\Sigma\setminus\gamma$ featuring in inequalities \eqref{eq:Panov-def-extended},\eqref{eq:p-function} become the values of $p_{\tilde u}$ on $\tilde \Sigma$ in inequalities \eqref{eq:Panov-def-extended-bis},\eqref{eq:p-function-bis} and vice versa.

Second, it is easily seen that
the entropy formulations \eqref{eq:Panov-def-extended},\eqref{eq:p-function} on $\Omega=(\R^+\times\R^d)\setminus\gamma$ and on $\R^+\times \R^d$ are in fact equivalent.
Indeed, sets in $\R^{d+1}$ of codimension at least $2$ have zero $\mathcal H^{d}$ Hausdorff measure and zero $W^{1,1}$ capacity.
 In particular, the measure in the right-hand side of \eqref{eq:Panov-def-extended} does not charge $\gamma$; the space $W^{1,1}_c((\R^+\times\R^d)\setminus\gamma)$ of Lipschitz continuous compactly supported functions, extended by $0$ on $\gamma$, is dense in $W^{1,1}_c(\R^+\times\R^d)$; and passage to the limit in
 \eqref{eq:Panov-def-extended} with a sequence of $W^{1,1}_c((\R^+\times\R^d)\setminus\gamma)$  test functions
 converging in $W^{1,1}(\R^+\times\R^d)$ leads to the same inequality with the limit test function.
Therefore, considering $W^{1,1}_c(\Omega)$ test functions and  considering
 $\mathcal D(\R^+\times\R^d)$ test functions in the context of integral inequalities
\eqref{eq:Panov-def-extended} is equivalent.

Combining the first and the second observation, we justify our claim. \end{proof}


Due to Proposition~\ref{prop:equiv-changeofvar}, existence of
admissible solutions for \eqref{eq:cons-law} with almost rectifiable
jump manifolds $\Sigma$ is reduced to existence of solutions, in the
sense of Definition~\ref{def-mman}, to analogous conservation law
set up in $\tilde \Omega$ with coefficients that may become singular
in a neighborhood of $\partial\tilde \Omega$; the crucial advantage of
\eqref{general_app4} is that it admits only flat jump manifolds.
Then, existence of admissible solutions to \eqref{general_app4} is
ensured by a combination of vanishing viscosity and truncation of $\tilde\Omega$.
First, one can approximate $\tilde\Omega$ from inside
 by a sequence of open domains $\tilde \omega_h$; their convergence to $\tilde\Omega$ is understood in the sense that
 $$
 \begin{array}{l}
 \text{$\cup_{\eps>0}\tilde\omega_\eps=\tilde\Omega$, $\cup_{\eps>0}\pa\tilde\omega_\eps\supset \tilde\Omega_0$,}\\
  \text{and $(\tilde\omega_\eps)_\eps$ are embedded in the sense that $\tilde\omega_{\eps_0}\subset\tilde\omega_{\eps}$ for $0<\eps<\eps_0$.}
 \end{array}
 $$
%
Then $\tilde \mff$ is globally Lipschitz continuous on
$\tilde\omega_\eps\times[a,b]$, $(t,x,k)\mapsto\tilde S(t,x,k)$ is
uniformly Lipschitz continuous in $k$ on $\tilde\omega_\eps\times[a,b]$.
We can choose $\tilde\omega_\eps$
with the property that the boundary of $\pa\tilde\omega_\eps$ is
piecewise parallel, piecewise orthogonal to the time axis. This can be achieved by covering compact subsets
of $\tilde \Omega\cup\tilde\Omega_0$ by a finite number of
cylinders contained in $\tilde\Omega\cup\tilde\Omega_0$. We will
denote by $\pa^{par}\tilde \omega_\eps$ the parabolic boundary of $\tilde
\omega_\eps$, i.e., the part of $\pa\tilde\omega_\eps$ where the exterior unit
normal vector is either orthogonal to the time direction (the
lateral boundary), or it points in the direction of increasing time
(this is the union of lower boundaries of $\tilde\omega_\eps$, that
contains $\tilde\Omega_0\cap\pa\tilde\omega_\eps$; but in general, it can
also contain pieces of hyperplanes $\{(t,\my)\,|\,t=const\}$). Then,
in the place of the initial-value viscosity regularized problem for the equation
\begin{equation}
\label{general_app4-epsilon}
 \pa_t \tilde{u} +\text{div}_{\my}
\tilde \mff (t,\my,\tilde{u})=\tilde S(t,\my,\tilde{u})+\eps \Delta_y \tilde{u},
\end{equation}
we consider the more general problem with the condition $u=a$
prescribed on the part
$\pa^{par}\tilde\omega_\eps\setminus\tilde\Omega_0$ of the parabolic
boundary:
 \begin{equation}\label{eq:parabolic-boundary-BC}
 \tilde  u|_{\pa^{par}\tilde\omega_\eps\setminus\tilde \Omega_0}=a
 \end{equation}
 and initial condition on the lower boundary:
 \begin{equation}\label{eq:lower-IC}
    \tilde u(0,\cdot)|_{\pa^{par}\tilde\omega_\eps\cap\tilde\Omega_0}=\tilde u_0\in [a,b].
 \end{equation}
 While this is not a
classical initial-boundary value problem because $\tilde\omega_\eps$ is
not cylindrical, solutions can be constructed in the same way as for
the classical setting of cylindrical domain,  (cf. the proof of Theorem~\ref{thm5}). To account for the specificity of the domain geometry, one can combine a Galerkin
semi-discretization scheme with a restarting procedure that accounts
for a finite number of changes of geometry of the sections $\tilde\omega_\eps^t:=\{\my\in\R^d\,|\,(t,\my)\in\tilde\omega_\eps\}$. We denote by $\tilde u_\eps$ the so constructed solution of \eqref{general_app4-epsilon},\eqref{eq:parabolic-boundary-BC},\eqref{eq:lower-IC}.

\smallskip
 The so constructed solutions obey the comparison
property, which is also proved in the same way as for the classical
initial-boundary value problems, cf. the proof of Lemma~\ref{l1}. Recall that due
to the choice of $[a,b]$-valued boundary and initial conditions,
 property \eqref{eq:a,b=sub,supersolution} is still valid for the
generalized initial-boundary value problem
\eqref{general_app4},\eqref{eq:parabolic-boundary-BC},\eqref{eq:lower-IC} set in
$\tilde\omega_\eps$. This guarantees that $\tilde u^\eps$ is $[a,b]$-valued.
Due to  the genuine nonlinearity property \eqref{gnl} of $\tilde \mff$
 and the embedding property of $(\tilde\omega_\eps)_\eps$, for all fixed $\eps_0>0$
$(\tilde u_\eps)_\eps$ has an accumulation point $u$ defined in $\tilde \omega_{\eps_0}$. Then, the diagonal extraction procedure and the convergence of $\tilde\omega_\eps$ to $\tilde\Omega$ in the sense explained above permit us to define $u$ globally on $\tilde\Omega$.

%

Therefore,  arguments developed in the proof of Theorem~\ref{thm5} and the related lemmas
prove that $\tilde u=\lim_{\eps\to 0} \tilde u_\eps$ verifies~\eqref{eq:Panov-def-extended-bis} with test functions
in $\mathcal D'(\tilde\omega_{\eps_0})$ and some values $p_{\tilde u}$ as in \eqref{eq:p-function-bis}.
It is easily seen that the initial data pass to the limit, i.e., $\tilde u(0,\cdot)=\tilde u_0$.
Due to convergence of $\tilde\omega_{\eps_0}$ to $\tilde\Omega$ as $\eps_0\to 0$, this ensures that
$u$ is a solution to \eqref{general_app4} in the sense of Definition~\ref{def-mman} with initial datum $\tilde u_0$.

Observe that the above approach gives not only an existence proof but also a uniqueness
proof, due to the equivalence stated in Proposition~\ref{prop:equiv-changeofvar} and to an elementary extension of
the calculations of \S~\ref{ssec:FLAT-direct-uniq} to the case of several flat interfaces (to be precise,
for the uniqueness of solutions in the sense of Definition~\ref{def-mman}, a uniform Lipschitz condition
of dependence on $u$ of $\tilde\mff$,$\tilde S$ should be imposed; it is verified in the examples we provide below).
Yet the direct uniqueness proof for solutions of \eqref{eq:cons-law} can be conducted without making appeal to a global change of variables, see \S~\ref{ssec:CURVED-prelim} and \cite{kar3,AKRnhm}. Therefore, we can summarize our result as follows.
\begin{theorem}\label{th:well-posedness-general}
  Consider equation~\eqref{eq:cons-law} with initial datum $u_0$. We impose:\\[3pt]
  $\bullet$ the confinement assumption \eqref{eq:confinement} on the flux $\mff$, source $S$ and initial datum $u_0$;\\
  $\bullet$ the genuine nonlinearity assumption~\eqref{gnl} of the flux $\mff$;\\
  $\bullet$ the almost-rectifiability assumption (Definition~\ref{def:almost-rec}) on the set $\Sigma$ of jump discontinuities in $(t,\mx)$ of the flux $\mff$;\\
  $\bullet$ on every connected component of $(R^+\times\R^d)\setminus\Sigma$, the assumptions of Lipschitz continuity of the flux $\mff$ with respect to all variables $(t,\mx,k)$ and of $\pa_u \mff$ with respect to $(t,\mx)$;\\
 $\bullet$  the assumption of global Lipschitz continuity of $k\mapsto S(t,\mx,k)$ for the source term.\\[3pt]
  Then there exists a unique weak solution to equation~\eqref{eq:cons-law} that is admissible in the sense of  \eqref{eq:Panov-def-extended},\eqref{eq:p-function} and takes the initial value $u_0$.
\end{theorem}

Now, we will provide two examples explaining how to check the almost-rectifiability
assumptions introduced above. Along with the example of Figure~\ref{fig:CircleToSquare}, this gives an idea of how large is the area of applicability of Theorem~\ref{th:well-posedness-general}.

\subsection{Two examples of almost rectifiable jump manifolds}\label{ssec:CURVED-rectifiable-examples}

 Now, we will explicitly construct the almost-rectifying maps $\Phi$ for two basic examples
 that permit to overview the ideas helpful in construction of almost-rectifying maps of Definition~\ref{def:almost-rec}.

 First, we deal with a family of curved non-crossing surfaces (for the case $d=2$, with jump surfaces parallel to the time axis) assuming that they have a common, globally defined transversal direction.
 The crossing is permitted in the second example, where we focus on the case of two merging curves (for the case $d=1$, with $t$-dependent jump curves). Clearly, the techniques for these examples can be combined.
 For example, the case of jump across the hyper-cylinder $\Sigma=\{(t,\mx)\,|\, t\geq 0,\, |\mx|=1\}$ (which is a prototype of the situation not covered by the first example) can be treated by considering $\Sigma$ as the union of two hypersurfaces that merge along the lines $\{(t,0,\hat x_1)\,|\,t\geq 0,\,|\hat x_1|=1\}$, in the spirit of the second example.

\subsubsection{The case of several non-intersecting simple discontinuity manifolds.}\label{sssec:CURVED-SeveralManifolds}

For simplicity, we will work in two-dimensional space, where we select the preferred direction $x_1$; we have $\hat x_1=x_2$. Assume that the flux discontinuity corresponding to the $t$-independent flux $\mff$ from
\eqref{eq:cons-law} consists of several curves (that become hypersurfaces, if the time variable is also considered) prescribed via Lipschitz continuous functions $\phi_i:\R\to \R$, $i=1,2,\dots,n$:
$$
\Sigma=\R^+\times(\cup_{i=1}^n\sigma_i),\;\; \sigma_i=\{x\in \R: \; x_1=\phi_i(x_2) \},
$$
with a strict separation (non-intersection) condition $\phi_{i}-\phi_{i-1}\geq \delta>0$ for $i=1,\dots,n-1$.

The equation that we are considering is thus
(with the notation $\chi_\Omega$ for the characteristic function of the set $\Omega$)
\begin{align}
\label{non-int} &\pa_t u+\sum\limits_{i=1}^{n+1}\text{div}_{\mx}(F_i(t,\mx,u)\chi_{\omega_i}(\mx))
=0,\\
\label{genic} &u|_{t=0}=u_0, \ \ a \leq u_0 \leq b,
\end{align} where $F_i(t,\mx,u)=(F_{1,i}(t,\mx,u),F_{2,i}(t,\mx,u))$
are Lipschitz continuous
and $\omega_i$ is the open region between the curves $\sigma_{i-1}$ and
$\sigma_{i}$, $i=1,\dots,n+1$, with the convention $\sigma_0=\{(x_1,x_2)\,|\,x_1=-\infty \}$ and
$\sigma_{n+1}=\{(x_1,x_2)\,|\,x_1=+\infty \}$.
We can also add a source term here, but we will omit it because it plays a passive role in the construction
of the almost-rectifying map for $\Sigma$.

\begin{figure}[htp]
\begin{center}
  \includegraphics[width=4in]{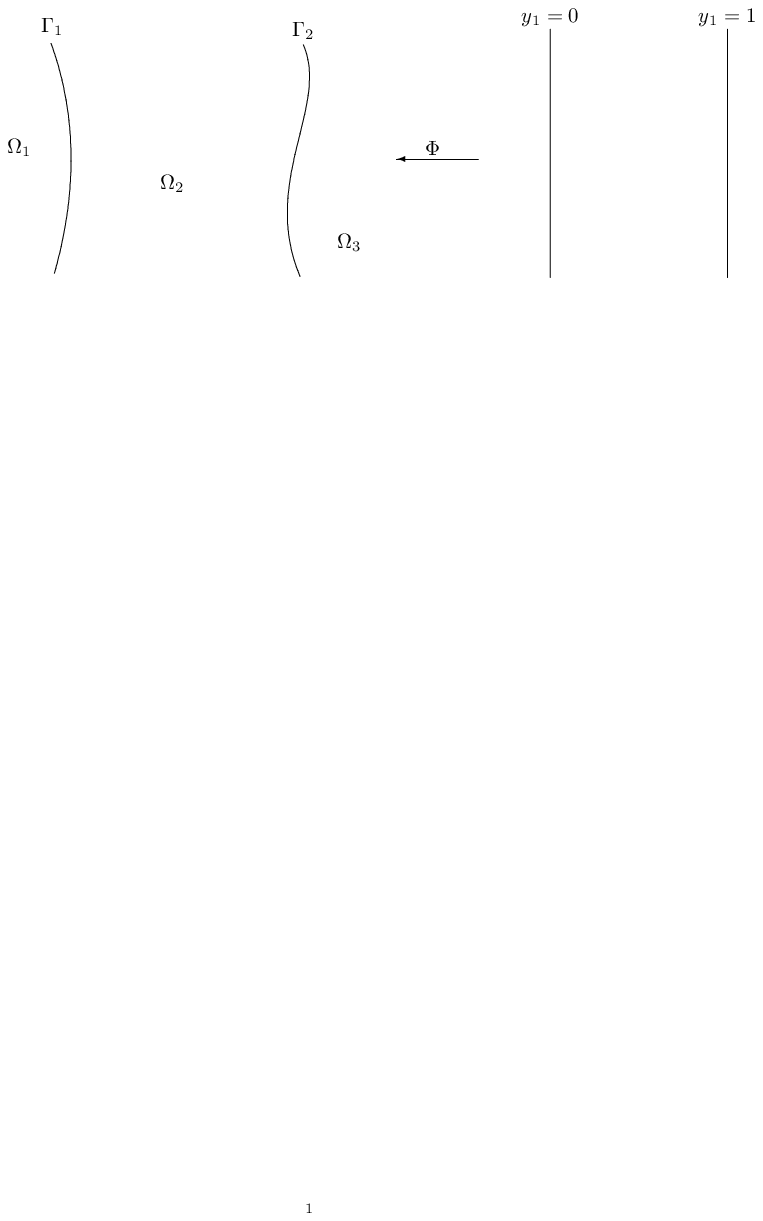}
  \caption{Case of several non-intersecting hypersurfaces in $\Sigma$}
  \label{fig:multi-surface}
  \end{center}
\end{figure}

 Now, the transformation $\Phi$ from Definition \ref{def:almost-rec} is constructed as follows
(see Figure~\ref{fig:multi-surface}).
In the domain $\omega_1=\{(x_1,x_2)\,|\,x_1\leq \phi_1(x_2) \}$ we set
$$
y_1=x_1-\phi_1(x_2), \ \ y_2=x_2.
$$
 In the domains $\omega_i=\{(x_1,x_2)\,|\, \phi_{i-1}(x_2)\leq x_1\leq \phi_i(x_2) \}$, $i\in \{2,\dots,n-1\}$, we set

$$
y_1=\frac{x_1-\phi_i(x_2)}{\phi_{i-1}(x_2)-\phi_i(x_2)}+i-2, \ \
y_2=x_2.
$$ In the domain $\omega_{n+1}\{(x_1,x_2)\,|\,\phi_n(x_2)\leq x_1 \}$, we take
$$
y_1=x_1-\phi_n(x_2)+(n-1), \ \ y_2=x_2.
$$
Since the curves are non-intersecting, the above defined transform
is bi-Lipschitz globally. Its inverse, given by is the formula
\begin{multline*}
x_1=(\phi_1(y_2)+y_1)H(-y_1)
\\
+ \sum_{i=1}^{n-1}\Bigl((\phi_{i+1}(y_2)-(y_1-i)(\phi_{i}(y_2)-\phi_{i+1}(y_2)) \Bigr) H(y_1-(i\!-\!1))H(i-y_1)\\
 + (\phi_n(y_2)+y_1-(n\!-\!1))H(y_1-(n\!-\!1))
\end{multline*}
with $H$ denoting the Heaviside function,
satisfies the properties required in Definition \ref{def:almost-rec} with $\gamma=\emptyset$, $\Gamma=\emptyset$.


\subsubsection{The case of two merging lines.}\label{sssec:CURVED-MergingLines}
\begin{figure}[htp]
\begin{center}
  \includegraphics[width=4in]{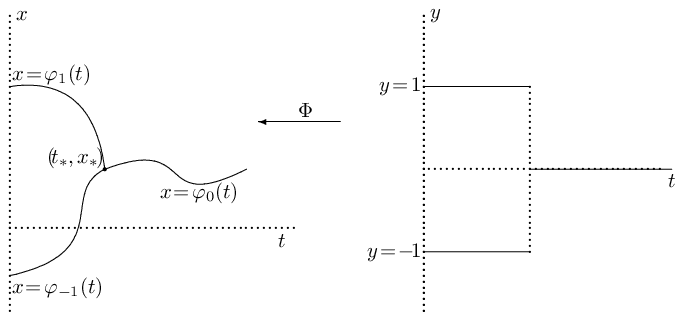}
  \caption{Case of two merging curves in $\Sigma$}
  \label{fig:merginglines}
  \end{center}
\end{figure}

Consider a scalar conservation law whose flux has discontinuities
disposed along two curves $\{(t,\phi_{-1}(t))\,|\, t\in[0,t_*]\}$ and $\{(t,\phi_1(t))\,|\, t\in[0,t_*]\}$
merging into one curve $\{(t,\phi_0(t))\,|\, t\in[t_*,+\infty)\}$ at the point $(t_*,x_*)$ (see Figure~\ref{fig:merginglines}).

Let $\phi_{-1},\phi_1,\phi_0$ be three functions of $t$ that
coincide at $t=t_*$, with the common value $x_*$. We assume that for
all $t<t_*$, $\phi_{-1}(t)<\phi_1(t)$. The functions  $\phi_{-1},\phi_1,\phi_0$ are assumed to be
Lipschitz continuous.


Then we define $\gamma:=\{(t_*,x_*) \}$, whose codimension with
respect to $\R^{d+1}$, $d=1$, is equal to two. We construct the map
from $(t,y)\in R^+\times \R$ defined as follows. The straight rays
or segments composing the set $\Gamma:=\{(t_*,y)\,|\, y\in [-1,1]\}$
and the set
$$\tilde \Sigma:=\Bigl([0,t_*)\times \{-1\}\Bigr)\cup\Bigl([0,t_*)\times \{+1\}\Bigr)
\cup \Bigl((t_*,+\infty)\times \{0\}\Bigr)
$$
split $\R^+\times\R$ into five regions $A$, $B$, $B'$, $C$ and $C'$. The map $\Phi: (t,y)\mapsto (t,x)=(t,x(t,y))$ is defined per region (we consider all the regions as closed, because the values of $\Phi$ match on the boundaries between the regions).
\begin{itemize}

  \item[$A.$] In the region $t\leq t_*$, $-1\leq y \leq 1$,
  $$
  x(t,y)=\phi_{-1}(t)+ \frac{y+1}{2} (\phi_1(t)-\phi_{-1}(t)).
  $$
  The image of this region is the space between the curves
  $\{(t,\phi_{-1}(t))\,|\, t\leq t_*\}$ and $\{(t,\phi_{1}(t))\,|\, t\leq t_*\}$.

  \item[$B.$] In the region $t\leq t_*$, $y\leq -1$,
  $$
  x(t,y)=\phi_{-1}(t) + (1+|t-t_*|)(y+1).
  $$
  The image is the region $\{(t,x)\,|\,t\leq t_*,\, x\leq \phi_{-1}(t)\}$.

  \item[$B'.$] In the region $t\geq t_*$, $y\leq 0$,
  $$
  x(t,y)= \phi_0(t) + (y+1) H(-(y+1) ) + |t-t_*|y.
  $$
  The image is the region  $\{(t,x)\,|\,t\geq t_*,\, x\leq \phi_{0}(t)\}$.

  \item[$C.$] In the region $t\leq t_*$, $y\geq 1$,
  $$
  x(t,y)=\phi_{1}(t) + (1+|t-t_*|)(y-1).
  $$
  The image is the region $\{(t,x)\,|\,t\leq t_*,\, x\geq \phi_{1}(t)\}$.

  \item[$C'.$] In the region $t\geq t_*$, $y\geq 0$,
  $$
  x(t,y)= \phi_0(t) + (y-1)H(y-1) + |t-t_*|y.
  $$
  The image is the region  $\{(t,x)\,|\,t\geq t_*, x\geq \phi_{0}(t)\}$.
\end{itemize}

 This is an example of map constructed according to the recipe of Remark~\ref{rem:changeofvar}.
 It is Lipschitz continuous (observe that each of the regions $B'$, $C'$ is split into two subregions
 with different polynomials defining the map) matching on the boundaries between (sub)regions.
 Both the map $\Phi$ and its inverse are Lipschitz away from any neighborhood of
 $\gamma$ or $\Gamma$, respectively. The second derivatives whose local boundedness is required in Definition~\ref{def:almost-rec}(a)(3) are easily computed and controlled.

\subsection{An alternative method for construction of solutions}\label{ssec:CURVED-LocalDiffusion}

Let us briefly indicate another way of constructing solutions.
\begin{definition}\label{def:local-rect}
A union of hypersurfaces $\Sigma$ in $\R^+\times\R^d$ is called locally almost rectifiable
if up to a set of codimension greater than $1$, $\Sigma$ can be included into the locally finite
union of disjoint open sets $U_i$ such that for all $i$, $\Sigma\cap U_i$ is a portion of the graph of a Lipschitz function $\psi_i:\R^d\mapsto \R$.
\end{definition}
This property is easily checked in practice.
For example, the cylinder $\Sigma=\{(t,x_1,x_2)\,|\, t\geq 0,\, x_1^2+x_2^2=1\}$ can be covered, up to the union
of the four lines $\{(t,0,\pm 1\}\cup\{(t,\pm 1,0\}$ and the curve $\{(0,x_1,x_2)\,|\, x_1^2+x_2^2=1\}$,
by the union of the four regions
$$U_i=\{(t,x_1,x_2)\,|\, t>0,\, \theta_i x_1>0,\, \theta^i x_2>0\,\},\;\; \theta_i,\theta^i\in \{-1,1\}.$$

In the sequel, we consider fluxes which singularities are included in a locally almost rectifiable set $\Sigma$.
According to assumption~\eqref{eq:time-orthog-Sigma}, we exclude the possibility that the discontinuity hypersurface is orthogonal to the time direction. Therefore up to an $i$-dependent rotation of space coordinates in $U_i$ we can assume that
 $$\Sigma\cap U_i=\{(t,x_1,\hat x_1)\,|\,x_1=\psi_i(t,\hat x_1)\},$$
  with Lipschitz continuous $\psi_i$.
The fact that $U_i$ are disjoint permits to consider a family of nonnegative $C^\infty_c(U_i)$ functions $(\lambda^\eps_i)_\eps$ with supports of $\lambda^\eps_i,\lambda^\eps_j$ that are disjoint, for all $i\neq j$,
and such that for all $i$, $\lambda^\eps_i$ converges to $1$ a.e. on ${U_i}$ as $\eps\to 0$; moreover, we can assume that $\lambda^\eps_i\equiv 1$ inside $U_i$ except in the $\eps$-neighbourhood of $\partial U_i$.

Now, we sketch the following construction that can be seen as another ``adapted viscosity'' approximation of \eqref{eq:cons-law}.

\subsubsection{An adapted diffusion operator.}\label{sssec:CURVED-AdaptedDiffusion}
 For every $U_i$, one can make the change of coordinates $ x_1=\tilde x_1+\psi_i(t,\hat x_1)$ that rectifies $\Sigma\cap U_i$; by analogy with the construction of Definition~\ref{def:almost-rec}, let us write $(t,\mx)=\Phi_i(t,\my)$, with $U_i=\Phi_i(\tilde U_i)$, but keeping in mind that $\Phi_i$ has a much simpler structure than in \S~\ref{ssec:CURVED-AlmostRect}
 (it is bi-Lipschitz, of jacobian $1$, and the mixed derivatives featuring in (a)(2) of Definition~\ref{def:almost-rec} are equal to zero). In the sequel, whenever $u$ is considered at $(t,\mx)\in U_i$,
 we mean that $(t,\mx)=\Phi_i(t,\my)$ with $(t,\my)\in \tilde U_i$, and we write $\tilde u$ for the transformed function $\tilde u (t,\my)=u(t,\mx)$.

 Denote by $\tilde \lambda^\eps_i$ the function in the transformed variables such that $\tilde \lambda^\eps_i(t,\my)=\lambda^\eps_i(t,\mx)$.
  In the new variables $\my$ in $\tilde U_i$, consider the following degenerate elliptic operator in conservative form:
   $$
   \tilde A^\eps_i: \tilde u\;\mapsto\; \eps\;\sum_{k=1}^d \partial_{y_k} (\tilde\lambda^\eps_i \partial_{y_k} \tilde u);
   $$
 observe that, by the definition of $\lambda^\eps_i$, $\tilde A^\eps_i$ acts as $\eps$ times the Laplacian in new variables $\my$ in a large portion of $\tilde U_i$.
Let $A^\eps_i$ be the corresponding operator acting on functions $u=\tilde u\circ \Phi_i^{-1}$ defined on $U_i$, this operator is implicitly defined by the change of variables $\Phi_i$.
For a global definition of $A^\eps_i$, observe that values of $(A^\eps_i u)(\cdot)$ outside $U_i$ can be set to zero, because $\lambda^\eps_i$ are supported in $U_i$.
Therefore $A^\eps_i$ is a degenerate anisotropic heterogeneous diffusion operator (in general, in non-divergence form) with coefficients that are Lipschitz continuous.

Due to the fact that the supports of $(\lambda^\eps_i)_i$ are disjoint, the operators $A^\eps_i$, implicitly defined on $U_i$ in the original coordinates, can be just pieced together to yield a second-order degenerate diffusion operator $A^\eps=\sum_i A^\eps_i$. By construction, near $\Sigma$ this operator represents the homogeneous isotropic Laplacian diffusion in coordinates adapted to the geometry of $\Sigma$, except in a $\eps$-neighbourhood $\mathcal N^\eps$ of $\Sigma\cap\Bigl(\cup_i \partial U_i\Bigr)$ in which the diffusion can degenerate.

Then we claim that the original equation \eqref{eq:cons-law} regularized with the following sum of diffusion operators:
 $$u\mapsto \delta\Delta u + A^\eps u$$
 permits to construct a solution of \eqref{eq:cons-law} in the sense of \eqref{eq:Panov-def-extended}, by letting first $\delta\downarrow 0$ with $\eps>0$ fixed, then letting $\eps\downarrow 0$.

Indeed, to start with, given $\delta,\eps>0$, existence of $L^2(\R^+;H^1_{loc}(\R^d))$ solutions $u^{\delta,\eps}$ with prescribed initial datum $u_0$
is ensured by classical theories of uniformly parabolic equations in non-divergence form.
It remains to pass to the limit and to characterize $\lim_{\eps\to 0}\lim_{\delta\to 0} u^{\delta,\eps}$
by inequalities \eqref{eq:Panov-def-extended}, with singular values \eqref{eq:p-function}.

\subsubsection{Auxiliary problem with degenerate adapted diffusion.}\label{sssec:CURVED-AuxiliaryAdaptedDiffProblem}
 In the first step, the compactness argument based on \eqref{gnl}, as used in \S~\ref{ssec:FLAT-exist}, permits to obtain an accumulation point $u^\eps$ from $(u^{\eps,\delta})_\delta$, as $\delta\to 0$.
 Away from an $\eps$-neighbourhood $\mathcal N^\eps$ of $\Sigma\cap\Bigl(\cup_i \partial U_i\Bigr)$, the function $u^\eps$ is an entropy solution  to equation \eqref{eq:cons-law} regularized with degenerate diffusion operator $A^\eps$:
 \begin{equation}\label{eq:cons-law-Adelta}
\partial_t u + \text{div}_\mx \mff (t,\mx,u) = A^\eps u + S(t,\mx,u).
\end{equation}
 Observe that in this limit problem, non-degenerate diffusion adapted to the geometry of $\Sigma$
 persists in a vicinity of $\Sigma\setminus \mathcal N^\eps$.

To be precise, the limit $u^\eps$ verifies in $(\R^+\times\R^d)\setminus \mathcal N^\eps$
the entropy formulation of \eqref{eq:cons-law-Adelta} analogous to the one put forward in the work \cite{BendahmaneKarlsen};
 in particular, $u^\eps$ is a local Kruzhkov entropy solution of \eqref{eq:cons-law}
in the regions where the flux $\mff$ is Lipschitz continuous and $A^\eps$ degenerates.
 The uniqueness theory of \cite{BendahmaneKarlsen} only covers the case of degenerate anisotropic diffusion but
 at a price of lengthy technicalities, it can be generalized in a straightforward way to
the heterogeneous case (cf. \cite{Vallet} for analogous extension of the isotropic homogeneous theory of \cite{Carrillo} to the isotropic heterogeneous case). Therefore, the (possibly local) formulation of \cite{BendahmaneKarlsen} for equation \eqref{eq:cons-law-Adelta} leads to the (local) Kato inequality for every couple of (local) entropy solutions $u^\eps,\hat u^\eps$ of \eqref{eq:cons-law-Adelta}.

\subsubsection{Convergence of adapted diffusion approximations.}\label{sssec:CURVED-ConvergenceAdapted}
  Starting from this point, our analysis mimics the one of \S~\ref{ssec:FLAT-exist}, using a family of explicit local viscosity profiles for \eqref{eq:cons-law-Adelta}; these are easily obtained in transformed coordinates $(t,\my)$,
  considering separately each of the sets $\tilde U_i$.

\smallskip
Namely, in the second step we apply the analogous \eqref{gnl}-based compactness argument to $(u^\eps)_\eps$.
An accumulation point $u$, as $\eps\to 0$, is a local Kruzhkov entropy solution of \eqref{eq:cons-law}
away from $\Sigma$. Indeed, the contribution of $A^\eps$ to the entropy dissipation is nonnegative,
and its contribution to the entropy flux vanishes, as $\eps\to 0$, because
there holds a uniform estimate on $\sum_i \|\sqrt{\eps} \lambda^\eps_i \nabla_\my  \tilde u^\eps \|_{L^2}$.
Further, since $U_i$ cover $\Sigma$ up to a lower-dimensional set, in order to  see that $u$ is a solution of \eqref{eq:cons-law} in the sense \eqref{eq:Panov-def-extended},\eqref{eq:p-function}
it is enough to justify the entropy inequalities \eqref{eq:Panov-def-extended} with test functions
supported in $U_i$, for every $i$, and some $p_u$ satisfying \eqref{eq:p-function}. The simplicity of the change of coordinates defined on $U_i$ implies
the invariance principle analogous to the one shown in Proposition~\ref{prop:equiv-changeofvar}: $u$
verifies \eqref{eq:Panov-def-extended},\eqref{eq:p-function} with test functions supported in $U_i$ if and only if $\tilde u$ verifies
the analogous inequalities with test functions supported in $\tilde U_i$ and appropriate singular values $p_{\tilde u}$. In variables $(t,\my)$ of $\tilde U_i$,
the interface $\tilde\Sigma_i:=\phi_i^{-1}(\Sigma\cap U_i)$ is flat, orthogonal to the direction $\tilde x_1$; and the diffusion that appears in the equation
on $\tilde u^\eps$ is the isotropic Laplacian $\eps\Delta_{\my}$ except in the set $\Phi^{-1}(U_i\cap \mathcal N^\eps)$ that vanishes, as $\eps\to 0$. Therefore in a neighbourhood of every point of $\tilde \Sigma_i$ one can construct profiles $\tilde R^\eps$ that only depend
on $\tilde x_1/\eps$, as in the construction of~\S~\ref{ssec:FLAT-exist}.
The analogue of Lemma~\ref{l1} follows from the Kato inequality involving $\tilde u^\eps$ and $\tilde R^\eps$.
Then or all $i$, the arguments of \S~\ref{ssec:FLAT-exist} permit to characterize the one-sided traces of $\tilde u$ on $\tilde \Sigma_i$ and to define singular values $p_{\tilde u}$, $\mathcal H^d$-a.e. on $\tilde \Sigma_i$. This provides entropy inequalities \eqref{eq:Panov-def-extended},\eqref{eq:p-function} with support of the test functions restricted to  $U_i$, where
the corresponding singular values $p_u$ are obtained from $p_{\tilde u}$ with the help of the transformation $\phi_i$.
Combining those with classical Kruzhkov entropy inequalities in $(\R^+\times\R^d)\setminus \Sigma$,
we justify \eqref{eq:Panov-def-extended},\eqref{eq:p-function} globally.

\smallskip
This concludes the sketch of construction of a solution $u$ in the sense \eqref{eq:Panov-def-extended},\eqref{eq:p-function}
under the assumption of local rectifiability of $\Sigma$ up to a lower-dimensional subset (Definition~\ref{def:local-rect}). In this way, the result fully analogous to the one of Theorem~\ref{th:well-posedness-general} can be obtained for
a different (at least as large as the class of almost rectifiable $\Sigma$'s, see Definition~\ref{def:almost-rec}) class of jump singularities $\Sigma$
in the flux $\mff$.

\section{Conclusion}
We have presented a new definition of solution to some discontinuous-flux problems,
along with detailed uniqueness and existence proofs for the case of a flat interface (a flux discontinuity hypersurface being seen as an interface).
 These results can be put in close correspondence with the results of \cite{AKRnhm}, where a simple particular case has been considered and quite different (rather artificial) approach to solutions characterization has been pursued.
 We think that, although our proofs in this paper do not appear any simpler than the proofs of \cite{AKRnhm}, the new definition of solution may be particularly useful in engineering applications because it complies with the physical and numerical intuition.
Further, we demonstrate that our well-posedness arguments extend to a  general configuration in spatially inhomogeneous media
with a locally finite number of Lipschitz-regular interfaces, by providing two original constructions of adapted viscosity approximate solutions.

\bigskip
\noindent
{\bf Acknowledgement} During the preparation of the article, Darko Mitrovic was engaged as a part time
researcher at the University of Bergen in the frame of the project
"Mathematical and Numerical Modeling over Multiple Scales" of the
Research Council of Norway whose support we gratefully acknowledge.
Boris Andreianov is supported by the French ANR project CoToCoLa.
The authors are grateful to the Berlin Mathematical School and to Prof. Etienne Emmrich from
the Technical University of Berlin for the invitation
for a research stay during which the paper has been written.

\end{document}